\documentclass{siamltex}
\usepackage{subeqn,extra}

\newtheorem{assumption}{Assumption}[section]
\newtheorem{condition}{Condition}[section]

\newcommand{\tE}{\tilde{E}}
\newcommand{\hE}{\hat{E}}
\newcommand{\tpert}{t}
\newcommand{\wDz}{\widehat{\Delta z}}
\newcommand{\wDl}{\widehat{\Delta \lambda}}
\newcommand{\wDs}{\widehat{\Delta s}}
\newcommand{\cirE}{\breve{E}}
\newcommand{\bE}{\bar{E}}
\newcommand{\be}{\bar{r}}
\newcommand{\cire}{\breve{r}}
\newcommand{\SN}{S_{\cN}}
\newcommand{\SB}{S_{\cB}}
\newcommand{\LB}{\Lambda_{\cB}}
\newcommand{\LN}{\Lambda_{\cN}}
\newcommand{\nablag}{{\nabla g}}
\newcommand{\nablagB}{{\nabla g_{\cB}}}
\newcommand{\nablagN}{{\nabla g_{\cN}}}
\newcommand{\alfmax}{\alpha_{\rm max}}
\newcommand{\walfmax}{\hat{\alpha}_{\rm max}}

\title{Effects of finite-precision arithmetic on 
interior-point methods for nonlinear programming}

\author{Stephen J. Wright%
\thanks{Mathematics and Computer Science Division,
Argonne National Laboratory, 9700 South Cass Avenue, Argonne, Illinois
60439, U.S.A.  This work was supported by the Mathematical,
Information, and Computational Sciences Division subprogram of the
Office of Advanced Scientific Computing, U.S. Department of
Energy, under Contract W-31-109-Eng-38.}}

\begin{document}
\markboth{STEPHEN J. WRIGHT}{FINITE-PRECISION EFFECTS IN NONLINEAR PROGRAMMING}
\pagestyle{myheadings}

\maketitle

%

\begin{abstract}
  We show that the effects of finite-precision arithmetic in forming
  and solving the linear system that arises at each iteration of
  primal-dual interior-point algorithms for nonlinear programming are
  benign, provided that the iterates satisfy centrality and
  feasibility conditions of the type usually associated with
  path-following methods. When we replace the standard assumption that
  the active constraint gradients are independent by the weaker
  Mangasarian-Fromovitz constraint qualification, rapid convergence
  usually is attainable, even when cancellation and roundoff errors
  occur during the calculations.
  In deriving our main results, we prove a key technical result about
  the size of the exact primal-dual step. This result can be used to
  modify existing analysis of primal-dual interior-point methods for
  convex programming, making it possible to extend the superlinear
  local convergence results to the nonconvex case.
\end{abstract}

\begin{AMS}
90C33, 90C30, 49M45
\end{AMS}

\section{Introduction} \label{sec.intro}

We investigate the effects of finite-precision arithmetic on the
calculated steps of primal-dual interior-point (PDIP) algorithms for
the nonlinear programming problem
\beq \label{nlp}
\makebox{\rm NLP:} \hspace*{.7in} 
\min_z \, \phi(z) \gap \makebox{\rm subject to $g(z) \le 0$},
\eeq
where $\phi:\R^n \to \R$ and $g:\R^n \to \R^m$ are twice Lipschitz
continuously differentiable functions. Optimality conditions for this
problem can be derived from the Lagrangian function $\cL(z,\lambda)$,
which is defined as
\beq \label{lagrange}
\cL(z,\lambda) 
= \phi(z) + \sum_{i=1}^m \lambda_i g_i(z) 
= \phi(z) + \lambda^T g(z),
\eeq
where $\lambda \in \R^m$ is a vector of Lagrange multipliers.  When a
constraint qualification (discussed below) holds at the point $z^*$,
first-order necessary conditions for $z^*$ to be a solution of
\eqnok{nlp} are that there exists a vector of Lagrange multipliers
$\lambda^* \in \R^m$ such that the following conditions are satisfied
for $(z,\lambda) = (z^*,\lambda^*)$:
\beq \label{kkt}
\cL_z(z,\lambda) = \nabla \phi(z) + \nabla g(z) \lambda = 0, \sgap
g(z) \le 0, \sgap
\lambda \ge 0, \sgap
\lambda^T g(z) = 0,
\eeq
where 
\[
\nabla g(z) = \left[ \nabla g_1(z), \nabla g_2(z), \dots, 
\nabla g_m(z) \right].
\]
The conditions \eqnok{kkt} are the well-known Karush-Kuhn-Tucker (KKT)
conditions.  We use $\cS_{\lambda}$ to denote the set of vectors
$\lambda^*$ such that $(z^*,\lambda^*)$ satisfies \eqnok{kkt}.  The
primal-dual solution set is defined by
\beq \label{pdsol}
\cS = \{ z^* \} \times \cS_{\lambda}.
\eeq

This paper discusses local convergence of PDIP algorithms for
\eqnok{nlp}, assuming that the algorithm is implemented on a computer
that performs calculations according to the standard model of
floating-point arithmetic. Because of our focus on {\em local}
convergence properties, we assume throughout that the current iterate
$(z,\lambda)$ is close enough to the solution set $\cS$ that
superlinear convergence would occur if exact steps (uncorrupted by
finite precision) were taken. In the interests of generality, we
weaken an assumption that is often made in the analysis of algorithms
for \eqnok{nlp}, namely, that the gradients of the active constraints
are linearly independent at the solution. We replace this linear
independence constraint qualification (LICQ) with the weaker
Mangasarian-Fromovitz constraint qualification (MFCQ)~\cite{ManF67}.
MFCQ allows constraint gradients to become dependent at the solution,
so that the set $\cS_{\lambda}$ of optimal Lagrange multipliers is no
longer necessarily a singleton, though it remains bounded. We continue
to assume that a strict complementarity (SC) condition holds, that is,
\beq \label{sc}
g_i(z^*)=0 \; \Rightarrow \; \lambda^*_i>0, \sgap 
\makebox{\rm for some $\lambda^* \in \cS_{\lambda}$}.
\eeq
In the context of rapidly convergent algorithms, the SC condition
makes good sense. If SC fails to hold, superlinear convergence of
Newton-like algorithms does not occur, except for
specially modified algorithms such as those that
identify the active constraints explicitly (see Monteiro and
Wright~\cite{MonW93a} and El-Bakry, Tapia, and Zhang~\cite{ElBTZ96}).

The major conclusion of the paper is that the effects of roundoff
errors on the rapid local convergence of the algorithm are fairly
benign.  When a standard second-order condition is added to the
assumptions already mentioned, the steps produced under floating-point
arithmetic approach $\cS$ almost as effectively as do exact steps, as
long as the distance to the solution set remains significantly greater
than the unit roundoff $\bu$.
The latter condition is hardly restrictive, since the data errors made
in storing the problem in a digital computer mean that the solution
set is known only to within some multiple of $\bu$ in any case.

The conclusions about the effectiveness of the computed steps are not
obvious, because all three formulations of the linear system that must
be solved to compute the step at each iteration may become highly
ill conditioned near the solution.
Our analysis would be significantly simpler if we were to make the
LICQ assumption because, in this case, one formulation of the linear
equations remains well conditioned, and stability of the three
standard formulations can be proved by exploiting their relationship
to this system of equations.

This work is related to earlier work of the author on finite-precision
analysis of interior-point algorithms for linear complementarity
problems~\cite{Wri93c} and linear
programming~\cite{Wri94a_rev,Wri96b}. The existence of second-order
effects gives the analysis here a somewhat different flavor, however.
In addition, we go into more depth in checking that the computed
iterates can continue to satisfy the approximate centrality conditions
usually required in primal-dual algorithms, and in deriving
expressions for the rate at which the computed iterates approach the
solution set. Related work by Forsgren, Gill, and
Shinnerl~\cite{ForGS94} deals with one formulation of the step
equations for the nonlinear programming problem---the so-called
augmented form treated here in Section~\ref{sec:pdip.aug}---but makes
assumptions on the pivot sequence that do not always hold in practice.
M.~H.~Wright~\cite{MWri98a} recently presented an analysis of the
condensed form of the step equations discussed in
Section~\ref{sec:pdip.cond} under the assumption that LICQ holds, and
found that the computed steps were more accurate than would be
expected from a naive analysis.

For linear programming, the PDIP approach has emerged as the most
powerful of the interior-point approaches. The supporting theory is
strong, in terms of global and local convergence analysis and
complexity theory (see the bibliography of Wright~\cite{IPPD96}).
Implementations yield better results than pure-primal or
barrier-function approaches; see Andersen et al.~\cite{AndGMX96}.
Strong theory is also available for these algorithms when applied to
convex programming, in which $\phi(\cdot)$ and $g_i(\cdot)$,
$i=1,\dots,m$ are all convex functions; see, for example, Wright and
Ralph~\cite{SJW32} and Ralph and Wright~\cite{RalW96,RalW96b}. The
latter paper drops the LICQ assumption in favor of  MFCQ,
making the local theory stronger in one sense than the corresponding
local theory for the sequential quadratic programming (SQP) algorithm.
The use of MFCQ complicates the analysis considerably, however; under LICQ, the
implicit function theorem can be used to prove a key technical result
about the length of the step, while more complicated logic is needed
to derive this same result under MFCQ.

A significant by-product of the current paper is to prove the key
technical result about the length of the rapidly convergent step
(Corollary~\ref{co:exact}) under MFCQ and SC, even when the problem
\eqnok{nlp} is not convex. This allows the local convergence results
of Ralph and Wright~\cite{SJW32,RalW96,RalW96b} to be extended to
general nonconvex nonlinear problems.


The analysis of this paper could also be applied to the recently
proposed stabilized sequential quadratic programming (sSQP) algorithm
(see Wright~\cite{Wri98a} and Hager~\cite{Hag97a}), in which small
penalties on the change in the multiplier estimate $\lambda$ from one
iteration to the next ensure rapid convergence even when LICQ is
relaxed to MFCQ. A finite-precision analysis of the sSQP method
appears in \cite[Section~3.2]{Wri98a}, but only for the augmented form
of the step equations. Analysis quite similar to that of the current
paper could be applied to show that similar conclusions continue to
hold when a condensed form of the step equations is used instead. We
omit the details.

The remainder of this paper is structured in the following way.
Section~\ref{notation} contains notation, together with our basic
assumptions about \eqnok{nlp} and some relevant results from the
literature. Section~\ref{PDIP} discusses the primal-dual
interior-point framework, defining the general form of each iteration
and the step equations that must be solved at each iteration.
Subsection~\ref{subsec.distance} proves an important technical result
about the relationship between the distance of an interior-point
iterate to the solution set $\cS$ and a duality measure $\mu$.
Section~\ref{accuracy} describes perturbed variants of the linear
systems that are solved to obtain PDIP steps, and proves our key
results about the effect of the perturbations on the accuracy of the
steps.

Section~\ref{sec:pdip.cond} focuses on one form of the PDIP step
equations: the most compact form in which most of the computational
effort goes into factoring a symmetric positive definite matrix,
usually by a Cholesky procedure. We trace the effect on step accuracy
of errors in evaluation of the functions, formation of the system, and
the factorization/solution process. Further, we show the effects of
these inaccuracies on the distance that we can move along the steps
before the interiority condition is violated, and on various measures
of algorithmic progress. An analogous treatment of the augmented form
of the step equations appears in Section~\ref{sec:pdip.aug}. The
conclusions of this section depend on the actual algorithm used to
solve the augmented system---it is not sufficient to assume, as in
Section~\ref{sec:pdip.cond}, that any backward-stable procedure is
used to factor the matrix. (We note that similar results hold for the
full form of the step equations, but we omit the details of this case,
which can be found in the technical report \cite{Wri98b}.)  We
conclude with a numerical illustration of the main results in
Section~\ref{sec:numerics} and summarize the paper in
Section~\ref{sec:conclusions}.


\section{Notation, Assumptions, and Basic Results} \label{notation}

We use $\cB$ to denote the set of active indices at $z^*$, that is,
\beq \label{Bdef}
\cB = \{ i=1,2,\dots,m \, | \, g_i(z^*) = 0 \},
\eeq
whereas $\cN$ denotes its complement 
\beq \label{Ndef}
\cN = \{ 1,2,\dots,m \} \backslash \cB.
\eeq
The set $\cB_+ \subset \cB$  is defined as 
\beq \label{B+def}
\cB_+ = \{ i \in \cB \, | \, \lambda^*_i >0 \;\; \makebox{\rm for some
$\lambda^*$ satisfying (\protect\ref{kkt})} \}.
\eeq
The strict complementarity condition \eqnok{sc} is equivalent to 
\beq \label{sc1}
\cB_+ = \cB.
\eeq
We frequently make reference to submatrices and subvectors
corresponding to the index sets $\cB$ and $\cN$. For example, the
quantities $\lB$ and $g_{\cB} (z)$ are the vectors containing the
components $\lambda_i$ and $g_i(z)$, respectively, for $i \in \cB$,
while $\nablagB(z)$ is the matrix whose columns are $\nabla g_i(z)$,
$i \in \cB$.

The Mangasarian-Fromovitz constraint qualification (MFCQ) is satisfied
at $z^*$ if there is a vector $\bar{y} \in \R^n$ such that 
\beq \label{mfcq} \nablagB(z^*)^T \bar{y} <0.  \eeq
The following fundamental result about MFCQ  is due to Gauvin~\cite{Gau77}.
\begin{lemma} \label{lem:Gauvin}
Suppose that the first-order conditions \eqnok{kkt}
are satisfied at $z=z^*$. Then $\cS_{\lambda}$ is bounded if and only if
the MFCQ condition \eqnok{mfcq} is satisfied at $z^*$.
\end{lemma}

This result is crucial because it allows our (local) analysis to place
a uniform bound on all $\lambda$ in a neighborhood of the dual
solution set $\cS_{\lambda}$.

The second-order condition used in most of the remainder of the paper
is that there is a constant $\xi>0$ such that
\beq \label{2os.L}
w^T \cL_{zz} (z^*, \lambda^*) w \ge \xi \| w \|^2,
\eeq
for all $\lambda^* \in \cS_{\lambda}$  and all $w$ satisfying
\beq \label{2os.1}
\begin{array}{ll} 
\nablag_i(z^*)^T w =0, & \makebox{\rm for all $i \in \cB_+$}, \\
\nablag_i(z^*)^T w \le 0, & \makebox{\rm for all $i \in \cB \backslash \cB_+$}.
\end{array}
\eeq
When the SC condition \eqnok{sc} (alternatively, \eqnok{sc1}) is
satisfied, this direction set is simply $\Null \nablag_{\cB}(z^*)^T$.

A simple example that satisfies MFCQ but not LICQ at the solution, and
that satisfies the second-order conditions \eqnok{2os.L},
\eqnok{2os.1} and the SC condition is as follows:
\beq \label{numex}
\min_{z \in \R^2} \, z_1 \;\; \mbox{subject to} \;\; 
(z_1-1/3)^2 + z_2^2 \le 1/9, \;\;
(z_1-2/3)^2 + z_2^2 \le 4/9.
\eeq
The solution is $z^*=0$, and the optimal multiplier set is
\beq \label{numex.sl}
\cS_{\lambda} = \{ \lambda \ge 0\, | \, 2 \lambda_1 + 4 \lambda_2 = 3 \}.
\eeq
The gradients of the two constraints are the solution are $(-2/3,0)^T$
and $(-4/3,0)^T$, respectively. They are linearly dependent, but the
MFCQ condition \eqnok{mfcq} can be satisfied by choosing $\bar{y} =
(1,0)^T$.

We use $\bu$ to denote the unit roundoff, which we define as the
smallest number such that the following property holds: When $x$ and
$y$ are any two floating-point numbers, $\mbox{\rm op}$ denotes $+$,
$-$, $\times$, or $/$, and $\mbox{\it fl}(z)$ denotes the
floating-point approximation of a real number $z$, we have
\beq \label{def.bu}
\mbox{\it fl}(x \; \mbox{\rm op} \; y) = 
(x \; \mbox{\rm op} \; y)(1+\epsilon), \sgap
|\epsilon| \le \bu.
\eeq
Modest multiples of $\bu$ are denoted by $\dbu$.

We assume that the problem is scaled so that the values of $g$ and
$\phi$ and their first and second derivatives in the vicinity of the
solution set $\cS$, and the values $(z,\lambda)$ themselves, can all
be bounded by moderate quantities. When multiplied by $\bu$,
quantities of this type are absorbed into the notation $\dbu$ in the
analysis below.

Order notation $O(\cdot)$ and $\Theta(\cdot)$ is used as follows: If
$v$ (vector or scalar) and $\epsilon$ (nonnegative scalar) are two
quantities that share a dependence on other variables, we write $v =
O(\epsilon)$ if there is a moderate constant $\beta_1$ such that $\| v
\| \le \beta_1 \epsilon$ for all values of $\epsilon$ that are
interesting in the given context. (The ``interesting context''
frequently includes cases in which $\epsilon$ is either sufficiently
small or sufficiently large, but we often use $v = O(\mu)$ to indicate
that $\| v \| \le \beta_1 \mu$ for all sufficiently small $\mu$ that
satisfy $\mu \gg \bu$, for some $\beta_1$; see later discussion.) We
write $v = \Theta(\epsilon)$ if there are constants $\beta_1$ and
$\beta_0$ such that $\beta_0 \epsilon \le \| v \| \le \beta_1
\epsilon$ for all interesting values of $\epsilon$.  Similarly, we
write $v = O(1)$ if $\| v \| \le \beta_1$, and $v = \Theta(1)$ if
$\beta_0 \le \| v \| \le \beta_1$.


We use the notation $\delta(z,\lambda)$ to denote the distance from
$(z,\lambda)$ to the primal-dual solution set, that is,
\beq \label{def.delta}
\delta(z,\lambda) \defeq 
\min_{(z^*,\lambda^*) \in \cS} \; \| (z,\lambda)-(z^*,\lambda^*) \|.
\eeq
It is well known (see, for example, Theorem~A.1 of Wright~\cite{Wri97e})
that this distance can be estimated in terms of known quantities at
$(z,\lambda)$. We state this result formally as follows.
\begin{theorem}\label{th:mu}
Suppose that the first-order conditions \eqnok{kkt}, the MFCQ
condition \eqnok{mfcq} and the second-order conditions
\eqnok{2os.L}, \eqnok{2os.1} are satisfied at the solution $z^*$.
Then if $\lambda \ge 0$, we have 
\beq \label{Sdist}
\delta(z,\lambda) = \Theta  \left( \left\|
\bmat{c}  \cL_z(z,\lambda)  \\ \min(\lambda, -g(z)) \emat
\right\| \right).
\eeq
\end{theorem}

We write the singular value decomposition (SVD) of the matrix
$\nablagB(z^*)$ of first partial derivatives as follows:
\beq \label{gen.3}
\nablagB (z^*) = \bmat{cc} \hat{U} & \hat{V} \emat
\bmat{cc} \Sigma & 0 \\ 0 & 0 \emat 
\bmat{c} U^T \\ V^T \emat =
\hat{U} \Sigma U^T,
\eeq
where the matrices $\bmat{cc} \hat{U} & \hat{V} \emat$ and $\bmat{cc}
U & V \emat$ are orthogonal, and $\Sigma$ is a diagonal matrix with
positive diagonal elements.  

Note that the columns of $\hat{U}$ form a basis for the range space of
$\nablagB(z^*)$, while the columns of $\hat{V}$ form a basis for the null
space of $\nablagB(z^*)^T$. Similarly, the columns of $U$ form a basis for
the range space of $\nablagB(z^*)^T$, while the columns of $V$ form a
basis for the null space of $\nablagB(z^*)$.   These
four subspaces are key to our analysis, particularly the subspace
spanned by the columns of $V$.  For the computational methods used to
solve the primal-dual step equations discussed in this paper, the
computed step in the $\cB$-components of the multipliers---that is,
$\Delta \lambda_{\cB}$---has a larger error in the range space of $V$
than in the complementary subspace spanned by the columns of $U$.  The
errors in the computed primal step $\Delta z$, the computed step of
the $\cN$-components of the multipliers $\lambda_{\cN}$, and the
computed step in the dual slack variables (defined later) are
typically also less significant than the error in $V^T \Delta
\lambda_{\cB}$. We show, however, that the potentially large error in
$V^T \Delta \lambda_{\cB}$ does not affect the performance of
primal-dual algorithms that use these computed steps until $\mu$
becomes similar to $\bu^{1/2}$. 

When the stronger LICQ condition holds, the matrix $V$
is vacuous, and the SVD \eqnok{gen.3} reduces to $\nablagB (z^*) =
\hat{U} \Sigma U^T$. Much of the analysis in this paper would simplify
considerably under LICQ, in part because $V^T \Delta
\lambda_{\cB}$---the step component with the large error---is no
longer present.

We use $\sigma_{\min}( \cdot )$ to denote the smallest eigenvalue,
and ${\rm cond}( \cdot)$ to denote the condition number, as measured
by the Euclidean norm.

\section{Primal-Dual Interior-Point Methods} \label{PDIP}

\subsection{Centrality Conditions and Step Equations}

Primal-dual interior-point methods are constrained, modified Newton
methods applied to a particular form of the KKT conditions
\eqnok{kkt}. By introducing a vector $s \in \R^m$ of slacks for the
inequality constraint, we can rewrite the nonlinear program as 
\[
\min_{(z,s)} \, \phi(z) \sgap \makebox{\rm subject to $g(z) + s=0$, $s \ge 0$,}
\]
and the KKT conditions \eqnok{kkt} as
\beq \label{kkt.s}
\cL_z(z,\lambda)  = 0, \sgap
g(z)+s=0, \sgap
(\lambda,s) \ge 0, \sgap
\lambda^T s=0.
\eeq
Motivated by this form of the conditions, we define
the mapping $\cF(z,\lambda,s)$ by
\beq \label{def.F}
\cF(z,\lambda,s) \defeq \bmat{c} \nabla \phi(z) + \nablag(z) \lambda \\
g(z) + s \\ S \Lambda e \emat,
\eeq
where the diagonal matrices $S$ and $\Lambda$ are defined by
\[
S \defeq \mbox{\rm diag}(s_1,s_2,\dots,s_m), \sgap
\Lambda \defeq \mbox{\rm diag}(\lambda_1,\lambda_2,\dots, \lambda_m),
\]
and $e$ is defined as 
\beq \label{def.e}
e=(1,1,\dots,1)^T.
\eeq
The KKT conditions \eqnok{kkt.s} can now be stated equivalently as
\beq \label{pdip.form}
\cF(z,\lambda, s) = 0, \gap (s, \lambda) \ge 0.
\eeq

Primal-dual iterates $(z,\lambda,s)$ invariably
satisfy the strict bound $(s,\lambda) >0$, while they
approach satisfaction of the condition $\cF(\cdot)=0$ in the
limit. An important measure of progress is the
{\em duality measure} $\mu(\lambda,s)$, which is defined  by
\beq \label{def.mu}
\mu (\lambda,s) \defeq \lambda^T s / m.
\eeq
When $\mu$ is used without arguments, we assume that $\mu =
\mu(\lambda,s)$, where $(z,\lambda,s)$ is the current primal-dual
iterate.  We emphasize that $\mu$ is a function of $(z,\lambda,s)$,
rather than a target value explicitly chosen by the algorithm, as is
the case in some of the literature.

A typical step $(\Dz, \Dl, \Ds)$ of the primal-dual method satisfies
\beq \label{pdip.step}
\nabla \cF (z,\lambda,s) \bmat{c} \Dz \\ \Dl \\  \Ds \emat =
-\cF(z,\lambda,s) - \bmat{c}  0 \\ 0 \\ \tpert \emat,
\eeq
where $\tpert$ defines the deviation from a pure Newton step for $\cF$
(which is also known as a ``primal-dual affine-scaling'' step). The
vector $\tpert$ frequently contains a centering term $\sigma \mu e$,
where $\sigma$ is a centering parameter in the range $[0,1]$.  It
sometimes also contains higher-order information, such as the product
$\Delta \Lambda_{\rm aff} \Delta S_{\rm aff} e$, where $\Delta
\Lambda_{\rm aff}$ and $\Delta S_{\rm aff}$ are the diagonal matrices
constructed from the components of the pure Newton step (Mehrotra~\cite{Meh92a}). 
In any case,
the vector $\tpert$ usually goes to zero rapidly as the iterates
converge to a solution, so that the steps generated from
\eqnok{pdip.step} approach pure Newton steps, which in turn ensures
rapid local convergence.  Throughout this paper, we assume that
$\tpert$ satisfies the estimate
\beq \label{tpert.est}
\tpert = O(\mu^2).
\eeq
All our major results continue to hold, with slight modification, if
we replace \eqnok{tpert.est} by $\tpert = O(\mu^{\sigma})$, for some
$\sigma \in (1,2]$. Our essential point remains unchanged; the
theoretical superlinear convergence rate promised by this choice of
$\tpert$ is not seriously compromised by roundoff errors as long as
$\mu$ remains significantly larger than the unit roundoff $\bu$.  To
avoid notational clutter, however, we analyze only the case
\eqnok{tpert.est}.

Using the definition \eqnok{lagrange}, we can write the system
\eqnok{pdip.step} explicitly as follows:
\beq \label{pdip.orig}
\bmat{ccc} \cL_{zz} (z,\lambda)  & \nablag(z) & 0 \\
\nablag(z)^T & 0 & I \\ 0 & S & \Lambda \emat 
\bmat{c} \Dz \\ \Dl \\ \Ds \emat =
- \bmat{c} \cL_{z} (z, \lambda) \\ g(z) + s \\ 
S \Lambda e + \tpert \emat.
\eeq
Block eliminations can be  performed  on this system to
yield more compact formulations. By
eliminating $\Ds$, we obtain the {\em augmented system}
form, which is
\beq \label{pdip.aug}
\bmat{cc} \cL_{zz} (z,\lambda) & \nablag(z) \\ 
\nablag(z)^T & -\Lambda^{-1} S \emat 
\bmat{c} \Dz \\ \Dl \emat =
\bmat{c} -\cL_{z} (z,\lambda)  \\ -g(z) + \Lambda^{-1} \tpert \emat.
\eeq
By eliminating $\Dl$ from this system, we obtain a system that is
sometimes referred to as the {\em condensed} form (or in the case of
linear programming as the {\em normal equations} form), which is
\beqa
\label{pdip.cond}
&& \left[ \cL_{zz} (z,\lambda) + \nablag(z) \Lambda S^{-1} \nablag(z)^T 
\right] \Dz \\
\nonumber
& = & -\cL_z (z,\lambda) - \nablag(z) \Lambda S^{-1}
[g(z) - \Lambda^{-1} \tpert ].
\eeqa

We consider primal-dual methods in which each iterate $(z,\lambda,s)$
satisfies the following properties:
\begin{subequations} \label{pdip.central}
\beqa
\label{pdip.central.1}
& \| r_f(z,\lambda) \| \le C \mu , \sgap \makebox{\rm where} \; 
r_f(z,\lambda) \defeq \cL_z (z,\lambda), \\
\label{pdip.central.2}
& \| r_g (z,s) \| \le C \mu, \sgap \makebox{\rm where} \; 
r_g (z,s) \defeq g(z) + s, \\
\label{pdip.central.3}
& (\lambda,s) >0, \sgap \lambda_i s_i \ge \gamma \mu, 
\sgap \makebox{\rm for all $i=1,2,\dots,m$,}
\eeqa
\end{subequations}
for some constants $C>0$ and $\gamma \in (0,1)$, where $\mu$ is
defined as in \eqnok{def.mu}. (In much of the succeeding discussion,
we omit the arguments from the quantities $\mu$, $r_f$, and $r_g$
when they are evaluated at the current iterate $(z,\lambda,s)$.)
These conditions ensure that the pairwise products $s_i \lambda_i$,
$i=1,2,\dots,m$ are not too disparate and that the first two
components of $\cF$ in \eqnok{def.F} can be bounded in terms of the
third component. They are sometimes called the {\em centrality
  conditions} because they are motivated by the notion of a central
path and its neighborhoods.  Conditions of the type
\eqnok{pdip.central} are imposed in most path-following interior-point
methods for linear programming (see, for example, \cite{IPPD96}).  For
nonlinear convex programming, examples of methods that require these
conditions can be found in Ralph and
Wright~\cite{SJW32,RalW96,RalW96b}. In nonlinear programming, we
mention Gould et al.~\cite{GouOST00} (see Algorithm 4.1 and Figure 5.1)
and Byrd, Liu, and Nocedal~\cite{ByrLN98}. In the latter paper,
\eqnok{pdip.central.1} and \eqnok{pdip.central.2} are imposed
explicitly, while \eqnok{pdip.central.3} can be guaranteed by choosing
$\epsilon_{\mu} = (1-\gamma) \mu$. Even when the choice
$\epsilon_{\mu} = \mu$ is made, as in the bulk of the discussion in
\cite{ByrLN98}, their other conditions concerning positivity of
$(s,\lambda)$ can be expected to produce iterates that satisfy
\eqnok{pdip.central.3} in practice.


For points $(z,\lambda,s)$ that satisfy \eqnok{pdip.central}, we can
use $\mu$ to estimate the distance $\delta(z,\lambda)$ from
$(z,\lambda)$ to the solution set $\cS$ (see \eqnok{def.delta}). These
results, which are proved in the following subsection, can be
summarized briefly as follows. When the MFCQ condition \eqnok{mfcq} and
the second-order conditions \eqnok{2os.L}, \eqnok{2os.1} are
satisfied, we have that $\delta(z,\lambda) = O(\mu^{1/2})$. When the
strict complementarity assumption \eqnok{sc} is added, we obtain the
stronger estimate $\delta(z,\lambda) = O(\mu)$. We can use these
estimates to obtain bounds on the elements of the diagonal matrices
$S$, $\Lambda$, $S^{-1} \Lambda$, and $\Lambda^{-1} S$ in the systems
above; these bounds are the key to the error analysis of the remainder
of the paper.


\subsection{Using the Duality Measure to  Estimate  Distance to the Solution}
\label{subsec.distance}

The main result of this section, Theorem~\ref{th:est3}, shows that
under certain assumptions, the distance $\delta(z,\lambda)$ of a
primal-dual iterate $(z,\lambda,s)$ to the solution set $\cS$ can be
estimated by the duality measure $\mu$. We start with a technical
lemma that proves the weaker estimate $\delta(z,\lambda) =
O(\mu^{1/2})$. 
Note that this result does not assume that the SC condition \eqnok{sc} holds.
\begin{lemma} \label{lem:est1}
Suppose that $z^*$ is a solution of \eqnok{nlp} at which the MFCQ
condition \eqnok{mfcq} and the second-order conditions \eqnok{2os.L},
\eqnok{2os.1} are satisfied. Then for all $(z,\lambda)$ with $\lambda
\ge 0$ for which there is a vector $s$ such that $(z,\lambda,s)$ 
satisfies \eqnok{pdip.central}, we have that
\beq \label{eq:11}
\delta(z,\lambda) = O(\mu^{1/2}).
\eeq
\end{lemma}
\begin{proof}
We prove the result by showing that
\beq \label{eq:12}
\bmat{c}  \cL_z(z,\lambda)  \\ \min(\lambda, -g(z)) \emat
 = O(\mu^{1/2})
\eeq
and then applying Theorem~\ref{th:mu}. Since $\cL_z(z,\lambda) = r_f =
O(\mu)$, the first part of the vector satisfies the required
estimate. For the second part, we have from \eqnok{pdip.central.2}
that
\[
-g(z) = s-r_g = s+O(\mu),
\]
and hence that
\beq \label{eq:12a}
\min(-g_i(z),\lambda_i) = \min(s_i,\lambda_i) + O(\mu).
\eeq
Because of \eqnok{def.mu} and \eqnok{pdip.central.3}, we have that
$s_i \lambda_i \le m \mu$ and $(\lambda_i,s_i) >0$. It follows
immediately that $\min(\lambda_i,s_i) \le (m \mu)^{1/2}$ for
$i=1,2,\dots,m$. Hence, by substitution into \eqnok{eq:12a}, we obtain
\[
\min(-g_i(z),\lambda_i) \le (m \mu)^{1/2} + O(\mu) = O(\mu^{1/2}).
\]
We conclude that the second part of the vector in \eqnok{eq:12} is of
size $O(\mu^{1/2})$, so the proof is complete.
\end{proof}

The following examples show the upper bound of Lemma~\ref{lem:est1} is
indeed achieved and that it is not possible to obtain a lower bound
on $\delta(z,\lambda)$ as a strictly increasing nonnegative
function of $\mu$. To demonstrate the
first claim, consider the problem
\[
\min \, \half z^2 \gap \makebox{\rm subject to $-z \le 0$.}
\]
The point $(z,\lambda,s) = (\epsilon, \epsilon,\epsilon)$ satisfies
\[
\cL_z(z,\lambda) = 0, \sgap
g(z) + s = 0, \sgap
s\lambda = \epsilon^2, \sgap \mu = \epsilon^2,
\]
so that the conditions \eqnok{pdip.central} are satisfied. Clearly the
distance from the point $(z,\lambda)$ to the solution set $\cS =
(0,0)$ is $\sqrt{2} \epsilon = \sqrt{2} \mu^{1/2}$. For the second
claim, consider any nonlinear program such that $\cB =
\{1,2,\dots,m\}$ (that is, all constraints active) and strict
complementarity \eqnok{sc} holds at some multiplier $\lambda^*$. 
Then for appropriate choices of $\gamma$ and $C$, the point
\beq \label{case2.pt}
(z,\lambda,s) = (z^*, \lambda^*, (m \mu) / (e^T \lamdba^*) e)
\eeq
satisfies the definition \eqnok{def.mu} and the condition
\eqnok{pdip.central} for any $\mu > 0$. On the other hand, we have
$\delta (z,\lambda)= \delta(z^*,\lambda^*)=0$ by definition, so there
are no $\beta>0$ and $\sigma>0$ that yield a lower bound estimate of
the form $\delta(z,\lambda) \ge \beta \mu^{\sigma}$.

We now prove an extension of Lemma~5.1 of Ralph and
Wright~\cite{RalW96}, dropping the monotonicity assumption of this
earlier result.
\begin{lemma} \label{lem:est2}
Suppose that the conditions of Lemma~\ref{lem:est1} hold and in
addition that the SC condition \eqnok{sc} is satisfied. Then for all
$(z,\lambda,s)$ satisfying \eqnok{pdip.central}, we have that
\begin{subequations} \label{eq:14}
\beqa
\label{eq:14a}
i \in \cB & \Rightarrow & s_i = \Theta(\mu), \;\; \lambda_i = \Theta(1), \\
\label{eq:14b}
i \in \cN & \Rightarrow & s_i = \Theta(1), \;\; \lambda_i = \Theta(\mu).
\eeqa
\end{subequations}
\end{lemma}
\begin{proof}
By boundedness of $\cS$ (Lemma~\ref{lem:Gauvin}), we have 
for all $(z,\lambda,s)$ sufficiently close to $\cS$ that
\beq \label{eq:15}
\lambda_i = O(1), \gap s_i = -g_i(z) + (r_g)_i = O(1).
\eeq
Given $(z,\lambda,s)$ satisfying \eqnok{pdip.central}, let $P(\lambda)$ be
the projection of $\lambda$ onto the set $\cS_{\lambda}$, and let  
$\lambda^* \in \cS_{\lambda}$  be some strictly complementary optimal
multiplier (for which \eqnok{sc} is satisfied). From Lemma~\ref{lem:est1}
we obtain
\beq \label{eq:16}
\| z-z^* \| = O(\mu^{1/2}).
\eeq
Using this observation together with smoothness of $\phi(\cdot)$ and
$g(\cdot)$, we have for the gradient of $\cL$ that
\beqas
\lefteqn{\cL_z(z,\lambda) - \cL_z(z^*,\lambda^*)} \\
&=& \nabla \phi(z) - \nabla \phi(z^*) + \nablag(z) \lambda -  \nablag(z^*) \lambda^* \\
&=& O(\mu^{1/2}) + \nablag(z) [\lambda - P(\lambda)] + 
[\nablag(z) - \nablag(z^*)]  P(\lambda) + 
\nablag(z^*) [P(\lambda) - \lambda^*].
\eeqas
Since $P(\lambda)$ and $\lambda^*$ are both in $\cS_{\lambda}$, we
find from \eqnok{kkt} that the last term vanishes. From \eqnok{eq:16}
and $P(\lambda) = O(1)$, the second-to-last term has size
$O(\mu^{1/2})$.  For the remaining term, we have $\nablag(z) = O(1)$, and
$\| \lambda - P(\lambda) \| \le \delta(z,\lambda) = O(\mu^{1/2})$. By
assembling all these observations, and using $\cL_z (z^*,\lambda^*)=0$, we obtain
\beq \label{eq:17}
\cL_z(z,\lambda) = \cL_z(z,\lambda) - \cL_z(z^*,\lambda^*) = O(\mu^{1/2}).
\eeq
Using again that $\nabla g(z^*) [ P(\lambda) - \lambda^*] = 0$,
we have from \eqnok{eq:16} that
\beqa
\nonumber 
[P(\lambda) - \lambda^*]^T [g(z) - g(z^*)] &=&
[P(\lambda) - \lambda^*]^T [\nablag(z^*)^T  (z-z^*)  + O(\| z-z^*\|^2)] \\
\label{eq:18}
&=& O(\| z-z^*\|^2)  = O(\mu).
\eeqa
By gathering the estimates \eqnok{eq:11},  \eqnok{eq:16}, \eqnok{eq:17}, and
\eqnok{eq:18}, we obtain
\beqa
\nonumber
\lefteqn{
\bmat{c} z-z^* \\ \lambda - \lambda^* \emat^T
\bmat{c} \cL_z (z,\lambda) - \cL_z (z^*,\lambda^*) \\ -g(z) + g(z^*) \emat } \\
\nonumber
&=&
\bmat{c} z-z^* \\ \lambda - P(\lambda) \emat^T
\bmat{c} \cL_z (z,\lambda) - \cL_z (z^*,\lambda^*) \\ -g(z) + g(z^*) \emat \\
\nonumber
&& 
+[P(\lambda) - \lambda^*]^T [-g(z) + g(z^*)] \\
\label{eq:19}
&=& O(\delta(z,\lambda)) O(\mu^{1/2}) + O(\mu) = O(\mu).
\eeqa
By substituting from \eqnok{pdip.central} and using \eqnok{eq:19}, we
obtain
\[
\bmat{c} z-z^* \\ \lambda - \lambda^* \emat^T
\bmat{c} r_f \\ s-r_g - s^* \emat = 
\bmat{c} z-z^* \\ \lambda - \lambda^* \emat^T
\bmat{c} \cL_z (z,\lambda) - \cL_z (z^*,\lambda^*) \\ -g(z) + g(z^*) \emat 
= O(\mu),
\]
and therefore 
\[
(\lambda - \lambda^*)^T (s-s^*) = -(z-z^*)^T r_f + (\lambda - \lambda^*)^T r_g
+ O(\mu).
\]
By using the conditions \eqnok{pdip.central.1}, \eqnok{pdip.central.2}, and the
definition \eqnok{def.mu}, we obtain
\beqas
\lefteqn{
-\sum_{i=1}^m \lambda^*_i s_i - \sum_{i=1}^m \lambda_i s^*_i} \\
&  = &
-(\lambda^*)^Ts - \lambda^T s^* = 
- \lambda^Ts + O(\mu) + O(\|z-z^*\|\|r_f\|) + 
O(\| \lambda - \lambda^*\| \| r_g\|)  = O(\mu).
\eeqas
Since $(\lambda,s)>0$ and $(\lambda^*,s^*) \ge 0$, all terms
$\lambda^*_i s_i$ and $\lambda_i s^*_i$, $i=1,2,\dots,m$ are
nonnegative, so there is a constant $C_1>0$ such that
\[
0 \le \lambda^*_i s_i \le C_1 \mu, \sgap
0 \le \lambda_i s^*_i \le C_1 \mu, \sgap 
\makebox{\rm for all $i=1,2,\dots,m$.}
\]
For all $i \in \cB$, we have $\lambda^*_i>0$ by our strictly
complementary choice of $\lambda^*$, and so
\beq \label{eq:20}
0 < s_i \le \frac{C_1}{\lambda^*_i} \mu \le
 \frac{C_1}{\min_{i \in \cB} \lambda^*_i} \mu \defeq C_2 \mu.
\eeq
On the other hand, we have by boundedness of $\cS_{\lambda}$ 
and our assumption \eqnok{pdip.central.3} that
\beq \label{eq:21}
s_i \ge \frac{\gamma \mu}{\lambda_i} \ge \gammin \mu, \sgap
\makebox{\rm for all $i=1,2,\dots,m$},
\eeq
for some constant  $\gammin>0$.
By combining \eqnok{eq:20} and \eqnok{eq:21}, we have that
\[
s_i = \Theta (\mu), \sgap \makebox{\rm for all $i \in \cB$}.
\]
For the $\lambda_{\cB}$ component, we have that
\[
s_i \lambda_i \ge \gamma \mu \;\; \Rightarrow \;\; 
\lambda_i \ge \frac{\gamma \mu}{s_i}  \ge \frac{\gamma}{C_2},
\sgap \makebox{\rm for all $i \in \cB$}.
\]
Hence, by combining with \eqnok{eq:15}, we obtain that
\[
\lambda_i = \Theta (1), \sgap \makebox{\rm for all $i \in \cB$}.
\]
This completes the proof of \eqnok{eq:14a}. We omit the proof of
\eqnok{eq:14b}, which is similar.
\end{proof}

Next, we show that when the strict complementarity assumption is added
to the assumptions of Lemma~\ref{lem:est1}, the upper bound on the
distance to the solution set in 
\eqnok{eq:11} can actually be improved to $O(\mu)$.
\begin{theorem} \label{th:est3}
Suppose that $z^*$ is a solution of \eqnok{nlp} at which the MFCQ
condition \eqnok{mfcq}, the second-order conditions \eqnok{2os.L},
\eqnok{2os.1}, and the SC condition \eqnok{sc} are
satisfied. Then for all $(z,\lambda)$ with $\lambda \ge 0$ for which
there is a vector $s$ such that $(z,\lambda,s)$ satisfies
\eqnok{pdip.central}, we have that
\beq \label{eq:25}
\delta(z,\lambda) = O(\mu).
\eeq
\end{theorem}
\begin{proof}
From \eqnok{pdip.central.1}, we have directly that $r_f = O(\mu)$.
We have from \eqnok{pdip.central} and \eqnok{eq:14a} that
\[
g_i(z) = -s_i + (r_g)_i = O(\mu), \sgap \lambda_i=\Theta(1), \sgap
\lambda_i>0 \sgap \makebox{\rm for all $i \in \cB$},
\]
so that
\beq \label{eq:27}
\min ( -g_i(z), \lambda_i) = -g_i(z) = O(\mu), 
\sgap \makebox{\rm for all $i \in \cB$,}
\eeq
whenever $\mu$ is sufficiently small.  For the remaining components,
we have
\beq \label{eq:28}
\min(-g_i(z) , \lambda_i) = \lambda_i = O(\mu), 
\sgap \makebox{\rm for all $i \in \cN$.}
\eeq
By substituting \eqnok{pdip.central.1}, \eqnok{eq:27}, and
\eqnok{eq:28} into \eqnok{Sdist}, we obtain the result.
\end{proof}


Similar conclusions to Lemma~\ref{lem:est2} and Theorem~\ref{th:est3}
can be reached in the case of linear programming algorithms.  The
second-order conditions \eqnok{2os.L}, \eqnok{2os.1} are not relevant
for this class of problems, and the SC assumption \eqnok{sc} holds for
every linear program that has a solution.

%

\section{Accuracy of PDIP Steps: General Results} \label{accuracy}

By partitioning the constraint index set $\{ 1,2,\dots,m\}$ into
active indices $\cB$ and inactive indices $\cN$, we can express the
system \eqnok{pdip.aug} without loss of generality as follows:
\beq \label{gen:aug}
\bmat{ccc} \cL_{zz} (z,\lambda) & \nablagB (z) & \nablagN(z) \\ 
\nablagB(z)^T & -D_{\cB} & 0 \\
\nablagN(z)^T & 0 & -D_{\cN} \\
\emat 
\bmat{c} \Dz \\ \DlB \\ \DlN \emat =
\bmat{c} -\cL_{z} (z,\lambda)  \\
-g_{\cB}(z) + \LB^{-1} \tpert_{\cB} \\
-g_{\cN}(z) + \LN^{-1} \tpert_{\cN}
\emat,
\eeq
where $D_{\cB}$ and $D_{\cN}$ are positive diagonal matrices defined by
\beq \label{def.DLS}
D_{\cB} = \LB^{-1} \SB , \gap D_{\cN} = \LN^{-1} \SN.
\eeq
When the SC condition \eqnok{sc} is satisfied, we have from
Lemma~\ref{lem:est2} that the diagonals of $D_{\cB}$ have size $\Theta(\mu)$
while those of $D_{\cN}$ have size $\Theta(\mu^{-1})$.
%
%
By eliminating $\DlN$ from \eqnok{gen:aug}, we obtain the following
intermediate stage between \eqnok{pdip.aug} and \eqnok{pdip.cond}:
\beqa \label{gen:augB}
\lefteqn{
\bmat{cc} H(z,\lambda) & \nablagB(z)  \\ 
\nablagB(z)^T & -D_{\cB} 
\emat 
\bmat{c} \Dz \\ \DlB \emat } \\
\nonumber
&=&
\bmat{c} -\cL_{z} (z,\lambda) - \nablagN(z) D_{\cN}^{-1}
[ g_{\cN}(z) -  \LN^{-1} \tpert_{\cN}] \\
-g_{\cB}(z) + \LB^{-1}   \tpert_{\cB}
\emat,
\eeqa
where we have defined
\beq \label{Hdef}
H(z,\lambda) \defeq \cL_{zz} (z,\lambda) + \nablagN(z) D_{\cN}^{-1} 
\nablagN(z)^T.
\eeq

In this section, we start by proving a key result about the solutions
of perturbed forms of the system \eqnok{gen:augB}. Subsequently, we
use this result as the foundation for proving results about the three
alternative formulations \eqnok{pdip.orig}, \eqnok{pdip.aug}, and
\eqnok{pdip.cond} of the PDIP step equations. The principal reason for
our focus on \eqnok{gen:augB} is that the proof of the main result can
be derived from fairly standard linear algebra arguments.
Gould~\cite[Section~6]{Gou86} obtains a system similar to
\eqnok{gen:augB} for the Newton equations for the primal log-barrier
function, and notes that the matrix approaches a nonsingular limit
when certain optimality conditions, including LICQ, are satisfied.
Because we replace LICQ by MFCQ, the matrix in \eqnok{gen:augB} may
approach a singular limit in our case.

We note that the form \eqnok{gen:augB} is also relevant to the
stabilized sequential quadratic programming (sSQP) method of
Wright~\cite{Wri98a} and Hager~\cite{Hag97a}; that is, slight
modifications to the results of this paper can be used to show that
the condensed and augmented formulations of the step equations for the
sSQP algorithm yield good steps even in the presence of roundoff errors and
cancellation. We omit further details in this paper.

Errors in the step equations arise from cancellation and roundoff
errors in evaluating both the matrix and right-hand side and from
roundoff errors that arise in the factorization/solution process. We
discuss these sources of error further and quantify them in the next
section.  In this section, we consider the following perturbed version
of \eqnok{gen:augB}:
\beq \label{gen.1}
\bmat{cc} H (z,\lambda) + \tE_{11} & \nablagB(z) + \tE_{12}  \\ 
\nablagB(z)^T + \tE_{21} & -D_{\cB} + \tE_{22}
\emat 
\bmat{c} w \\ y \emat
=
\bmat{c} r_1 \\ \nablagB(z^*)^T r_3 + r_4 \emat.
\eeq
Here, $\tE$ is the perturbation matrix (appropriately partitioned and
not assumed to be symmetric) and $r_1$, $r_3$, and $r_4$ represent
components of a general right-hand side. Note the partitioning of the
second right-hand side component into a component $\nablagB(z^*)^T r_3$ in
the range space of $\nablagB(z^*)^T$ and a remainder term $r_4$. When LICQ
is satisfied, the range space of $\nablagB(z^*)^T$ spans the full space,
so we can choose $r_4$ to be zero. Under MFCQ, however, we have in
general that $r_4$ must be nonzero. The main interest of the results
below is in isolating the component of the solution of \eqnok{gen.1}
that is sensitive to $r_4$.

To make the results applicable to a wider class of linear systems, we
do not impose the assumptions that were needed in the preceding section
to ensure that the matrices $D_{\cB}$ and $D_{\cN}$ defined by
\eqnok{def.DLS} have diagonals of the appropriate size. Instead, we
{\em assume} that the diagonals have the given size, and derive the
application to the linear systems of interest (those that arise in
primal-dual interior-point methods) as a special case.

Our results in this and later sections make extensive use of the SVD
\eqnok{gen.3} of $\nablagB (z^*)$. They also make assumptions about
the size of the smallest singular value of this matrix, and about the
size of the smallest eigenvalue of $\hat{V}^T \cL_{zz}(z^*,\lambda^*)
\hat{V}$, the two-sided projection of the Lagrangian Hessian onto the
active constraint manifold.

\begin{theorem} \label{th:gen1}
Let $(z,\lambda)$ be an approximate primal-dual solution of
\eqnok{nlp} with $\delta(z,\lambda) = O(\mu)$, and suppose the
diagonal matrices $D_{\cB}$ and $D_{\cN}^{-1}$ defined by
\eqnok{def.DLS} have all their diagonal elements of size
$\Theta(\mu)$.
Suppose that the perturbation submatrices in \eqnok{gen.1} satisfy
\beq \label{gen.2}
\tE_{11} = \dbu / \mu + O(\mu), \sgap 
\tE_{21}, \tE_{12}, \tE_{22} = \dbu,
\eeq
and that the following conditions hold for some $\beta>0$:
\begin{subequations} \label{gen.0}
\beqa \label{gen.0.1}
&  \bu / \mu \ll 1, \sgap  \bu \ll 1, \\
\label{gen.0.3}
& \sigma_{\min} (\Sigma) \ge \beta \max( \mu^{1/3}, \bu/\mu), \\
\label{gen.0.2}
& \sigma_{\min} (\hat{V}^T \cL_{zz}(z^*,\lambda^*) \hat{V}) 
\ge \beta \max(\mu^{1/3},\bu / \mu), 
\sgap \makebox{\rm for all $\lambda^* \in \cS_{\lambda}$}.
\eeqa
\end{subequations}
Then if $\beta$ is sufficiently large (in a sense to be specified in the
proof),
the step $(w,y)$ computed from \eqnok{gen.1} satisfies
\beqas
(U^Ty, \hat{V}^Tw, \hat{U}^T w) &=& O(\| r_1 \| + \| r_3 \| + \| r_4 \|), \\
V^Ty &=&  O(\| r_1 \| + \| r_3 \| + \| r_4\|/\mu).
\eeqas
\end{theorem}
\begin{proof}
If we define
\[
y_U = U^Ty, \sgap
y_V = V^Ty, \sgap
w_{\hat{U}} = \hat{U}^T w, \sgap
w_{\hat{V}} = \hat{V}^T w,
\]
we have
\[
y = Uy_U + V y_V, \gap
w = \hat{U} w_{\hat{U}} + \hat{V} w_{\hat{V}}.
\]
Using this notation, we can rewrite \eqnok{gen.1} as
\beqa 
\label{gen.5}
& \bmat{cccc} \hat{U}^T M_{11} \hat{U} & \hat{U}^T M_{11} \hat{V} &
\hat{U}^T M_{12} U & \hat{U}^T M_{12} V \\
\hat{V}^T M_{11} \hat{U} & \hat{V}^T M_{11} \hat{V} &
\hat{V}^T M_{12} U & \hat{V}^T M_{12} V \\
U^T M_{21} \hat{U} & U^T M_{21} \hat{V} & U^T M_{22} U & U^T M_{22} V \\
V^T M_{21} \hat{U} & V^T M_{21} \hat{V} & V^T M_{22} U & V^T M_{22} V 
\emat
\bmat{c} 
w_{\hat{U}} \\ w_{\hat{V}} \\y_U \\ y_V \emat \\
\nonumber
& \gap = \bmat{c} \hat{U}^T r_1 \\ \hat{V}^T r_1 \\ U^T \nablagB(z^*)^T r_3 +
U^T r_4 \\ V^T \nablagB(z^*)^T r_3 + V^T r_4 \emat,
\eeqa
where we have defined
\beqa
\label{Mdef}
& M_{11} = H (z,\lambda) + \tE_{11}, \sgap
M_{12} = \nablagB(z) + \tE_{12} , \\
\nonumber
& M_{21} = \nablagB(z)^T + \tE_{21}, \sgap
M_{22} =  -D_{\cB} + \tE_{22},
\eeqa
and $H(\cdot,\cdot)$ is defined in \eqnok{Hdef}. From \eqnok{gen.3}, we have
\[
V^T \nablagB(z^*)^T = 0, \sgap
U^T \nablagB(z^*)^T = \Sigma \hat{U}^T.
\]
The fact that $V^T$ annihilates $\nablagB(z^*)^T$ is crucial, because it
causes the term with $r_3$ to disappear from the last component of the
right-hand side of \eqnok{gen.5}, which becomes
\beq \label{gen.5a}
\bmat{c} \hat{U}^T r_1 \\ \hat{V}^T r_1 \\ \Sigma \hat{U}^T r_3 +
U^T r_4 \\ V^T r_4 \emat.
\eeq
From the definitions \eqnok{Mdef} and \eqnok{Hdef}, the perturbation
bound \eqnok{gen.2}, our assumptions that $D_{\cN}^{-1} = O(\mu)$ and
$\delta(z,\lambda) = O(\mu)$, compactness of $\cS$, and the fact that
$\cL_{zz}$ is Lipschitz continuous, we have that
\beq \label{M11.def}
M_{11} = \cL_{zz}(z^*,\lambda^*) + \dbu/\mu + O(\mu),
\eeq
for some $\lambda^* \in \cS_{\lambda}$.  Using these same facts, we
have likewise that
\[
M_{21} = \nablagB (z^*)^T + \dbu + O(\mu),
\]
so by substituting from \eqnok{gen.3}, we have that
\begin{subequations} \label{M21.defs}
\beqa
& U^T M_{21} \hat{U} = \Sigma + \dbu + O(\mu), \sgap
U^T M_{21} \hat{V} = \dbu + O(\mu), \\
& V^T M_{21} \hat{U} = \dbu + O(\mu), \sgap 
V^T M_{21} \hat{V} = \dbu + O(\mu).
\eeqa
\end{subequations}
Similarly, from the definition of $M_{12}$, we have
\begin{subequations} \label{M12.defs}
\beqa
& \hat{U}^T M_{12} U = \Sigma + \dbu + O(\mu), \sgap
\hat{U}^T M_{12} V = \dbu + O(\mu), \\
& \hat{V}^T M_{12} U = \dbu + O(\mu), \sgap
\hat{V}^T M_{12} V = \dbu + O(\mu).
\eeqa
\end{subequations}
For the $M_{22}$ block, we have from \eqnok{Mdef} and \eqnok{gen.2} that
\begin{subequations} \label{M22.defs}
\beqa
& U^T M_{22} U  = -U^T D_{\cB} U + \dbu = O(\mu) + \dbu, \\
& U^T M_{22} V  = O(\mu) + \dbu, \sgap
V^T M_{22} U = O(\mu) + \dbu, \\
& V^T M_{22} V = -V^T D_{\cB} V + \dbu =  \tilde{M}_{VV} + \dbu,
\eeqa
\end{subequations}
where $\tilde{M}_{VV} \defeq -V^T D_{\cB} V$ 
has all its singular values of size $\Theta(\mu)$, so that
\beq \label{mtilinv}
\tilde{M}_{VV}^{-1} = \Theta (\mu^{-1}).
\eeq
Using these estimates together with \eqnok{gen.5a}, 
we can rewrite \eqnok{gen.5} as
\beq \label{gen.10}
\left\{ \bmat{cc} Q & 0 \\ 0 & \tilde{M}_{VV} \emat  +
\bmat{cc} \hE_{11} & \hE_{12} \\ \hE_{21} & \hE_{22} \emat \right\}
\bmat{c} w_{\hat{U}} \\ w_{\hat{V}} \\y_U  \\ \hline y_V \emat  =
\bmat{c} \hat{U}^T r_1 \\ \hat{V}^T r_1 \\ \Sigma \hat{U}^T r_3 +
U^T r_4 \\ \hline V^T r_4 \emat,
\eeq
where
\beqa 
\label{gen.10b}
Q & = & \bmat{ccc} \hat{U}^T \cL_{zz} (z^*,\lambda^*) \hat{U} &
\hat{U}^T \cL_{zz} (z^*,\lambda^*) \hat{V} & \Sigma \\
\hat{V}^T \cL_{zz} (z^*,\lambda^*) \hat{U} &
\hat{V}^T \cL_{zz} (z^*,\lambda^*) \hat{V} & 0 \\
\Sigma & 0 & 0 \emat \\
\nonumber
&& +
\bmat{ccc}
\dbu/\mu + O(\mu) & \dbu/\mu + O(\mu) & \dbu + O(\mu)  \\
\dbu/\mu + O(\mu) & \dbu/\mu + O(\mu) & 0 \\
\dbu + O(\mu) & 0 & 0 
\emat \\
\label{gen.10a}
& \defeq &
\bmat{ccc} N_{UU} & N_{UV} & \bar{\Sigma}_1 \\
 N_{VU} & N_{VV} & 0 \\
\bar{\Sigma}_2 & 0 & 0 \emat,
\eeqa
while 
\beq \label{gen.11}
\hE_{11} = \bmat{ccc}
0 & 0 & 0 \\
0 & 0 & \dbu + O(\mu) \\
0 & \dbu + O(\mu) & \dbu + O(\mu) 
\emat,
\eeq
and
\beq \label{hatE}
\hE_{12}, \hE_{21} = \dbu + O(\mu) = O(\mu), \sgap
\hE_{22} = \dbu.
\eeq

For purposes of specifying the required conditions on $\beta$ in
\eqnok{gen.0.3} and \eqnok{gen.0.2}, we define $\kappa$ to be a
constant such that expressions of size $\dbu$ and $O(\mu)$ that
arise in the perturbation terms in the coefficient matrix in
\eqnok{gen.10} can be bounded by $\kappa \bu$ and $\kappa \mu$,
respectively. For example, we suppose that the perturbations in
$\bar{\Sigma}_1$, $\bar{\Sigma}_2$, and $N_{VV}$ can bounded as
follows:
\begin{subequations} \label{final1}
\beqa
\label{final1.1}
& \| \bar{\Sigma}_1 - \Sigma \| \le \kappa (\mu  + \bu), \sgap
\| \bar{\Sigma}_2 - \Sigma \| \le \kappa (\mu  + \bu), \\
& \| N_{VV} - \hat{V}^T \cL_{zz} (z^*,\lambda^*) \hat{V} \| \le 
\kappa (\bu / \mu + \mu),
\eeqa
\end{subequations}
and that 
\beq \label{final2}
\| \hat{E}_{11} \| \le \kappa (\bu + \mu), \;\;
\| \hat{E}_{12} \| \le \kappa (\bu + \mu), \;\;
\| \hat{E}_{21} \| \le \kappa (\bu + \mu), \;\;
\| \hat{E}_{22} \| \le \kappa \bu.
\eeq
From \eqnok{final1.1} and \eqnok{gen.0.3}, we have that
\[
\| \bar{\Sigma}_1 - \Sigma \| \le \kappa \max(\mu^{1/3}, \bu/\mu) \le
(\kappa / \beta) \sigma_{\rm min} (\Sigma) \le (\kappa / \beta) \| \Sigma \|.
\]
It is therefore easy to show that if $\beta$ can be chosen large
enough that $\beta > 2 \kappa$ (while still satisfying \eqnok{gen.0.3}
and \eqnok{gen.0.2}), then 
\beq \label{gen.60.a}
\| \bar{\Sigma}_1 \| \le 2 \| \Sigma \|, \sgap
\| \bar{\Sigma}_1^{-1} \| \le 2 \| \Sigma^{-1} \|.
\eeq
Similarly, we can show that
\beq
\label{gen.60.b}
\| \bar{\Sigma}_2 \| \le 2 \| \Sigma \|, \sgap
\| \bar{\Sigma}_2^{-1} \| \le 2 \| \Sigma^{-1} \|,
\eeq
\beq 
\label{gen.60.c}
\| N_{VV} \| \le 
    2 \| \hat{V}^T \cL_{zz} (z^*,\lambda^*) \hat{V} \|, \sgap
\| N_{VV}^{-1} \| \le 
    2 \| (\hat{V}^T \cL_{zz} (z^*,\lambda^*) \hat{V})^{-1} \|.
\eeq
Note, too, that because of Lipschitz continuity of $\cL_{zz}$ and
compactness of $\cS$, and the bounds \eqnok{gen.0.1}, the norms of
$N_{UU}$, $N_{UV}$, $N_{VU}$, $N_{VV}$, and $\Sigma$ are all $O(1)$.
Hence the matrix $Q$ is itself invertible, and we have
\beq \label{gen.12}
Q^{-1} = \bmat{ccc}
0 & 0 & \bar{\Sigma}_2^{-1} \\
0 & N_{VV}^{-1} & - N_{VV}^{-1}N_{VU} \bar{\Sigma}_2^{-1} \\
\bar{\Sigma}_1^{-1} & 
-\bar{\Sigma}_1^{-1} N_{UV} N_{VV}^{-1} & 
-\bar{\Sigma}_1^{-1} ( N_{UU} - N_{UV} N_{VV}^{-1} N_{VU}) \bar{\Sigma}_2^{-1} 
\emat.
\eeq
Noting that
\beq \label{gen.13}
(Q + \hE_{11})^{-1} = (I+Q^{-1} \hE_{11})^{-1} Q^{-1},
\eeq
we examine the size of $Q^{-1} \hE_{11}$. Note first from
\eqnok{gen.0.3} and \eqnok{gen.0.2} together with \eqnok{gen.60.a},
\eqnok{gen.60.b}, and \eqnok{gen.60.c} that
\begin{subequations} \label{gen.14}
\beqa \label{gen.14.a}
& \| \bar{\Sigma}_1^{-1} \| \le \frac{2}{\beta} ({\bu}/\mu)^{-1}, \sgap 
\| \bar{\Sigma}_2^{-1} \| \le \frac{2}{\beta} ({\bu}/\mu)^{-1}, \sgap 
\| N_{VV}^{-1} \| \le \frac{2}{\beta} ({\bu}/\mu)^{-1}, \\
\label{gen.14.b}
& \| \bar{\Sigma}_1^{-1} \| \le \frac{2}{\beta} \mu^{-1/3}, \sgap 
\| \bar{\Sigma}_2^{-1} \| \le \frac{2}{\beta} \mu^{-1/3}, \sgap 
\| N_{VV}^{-1} \| \le \frac{2}{\beta} \mu^{-1/3}.
\eeqa
\end{subequations}
By forming the product of \eqnok{gen.12} with \eqnok{gen.11} and using
the bounds in \eqnok{gen.14}, we can show that the norm of $Q^{-1}
\hE_{11}$ can be made less than $1/2$ provided that $\beta$ in 
\eqnok{gen.0.3}, \eqnok{gen.0.2} is 
sufficiently large. The $(3,3)$ block of $Q^{-1} \hE_{11}$, for
instance, has the form
\[
-\bar{\Sigma}_1^{-1} N_{UV} N_{VV}^{-1} (\dbu + O(\mu)) +
\bar{\Sigma}_1^{-1} (N_{UU} - N_{UV} N_{VV}^{-1} N_{VU}) \bar{\Sigma}_2^{-1} 
(\dbu + O(\mu)).
\]
Because of  \eqnok{final2}, its norm can be  bounded by a quantity of the form
\[
C \kappa \left( \| \bar{\Sigma}_1^{-1} \| \, \| N_{VV}^{-1} \|  + 
\| \bar{\Sigma}_1^{-1} \| \, \| \bar{\Sigma}_2^{-1} \| \, \| N_{VV}^{-1} \| +
\| \bar{\Sigma}_1^{-1} \| \, \| \bar{\Sigma}_2^{-1} \| \right) 
\left( (\bu/\mu) \mu + \mu \right),
\]
(for some $C$ that depends on $\| \cL_{zz}(z^*,\lambda^*) \|$), 
which in turn because of \eqnok{gen.14} is bounded by the following
quantity:
\[
8 C \kappa \left( \frac{1}{\beta^2} \mu^{2/3} + \frac{1}{\beta^3} \mu^{1/3} \right) +
8 C \kappa \left( \frac{1}{\beta^2} \mu^{1/3} + \frac{1}{\beta^3} \right).
\]
Provided that $\beta$ is large enough that this and the other blocks
of $Q^{-1} \hE_{11}$ can be bounded appropriately, we have that $\|
Q^{-1} \hE_{11} \| \le 1/2$, and therefore from \eqnok{gen.13} we have
\[
\| (Q + \hE_{11})^{-1}\|  = 2 \|  Q^{-1} \|.
\]
Our conclusion is that for $\beta$ satisfying the conditions outlined
in this paragraph, the inverse of the $(1,1)$ block of the matrix in
\eqnok{gen.10} can be bounded in terms of $\| Q^{-1} \|$, which
because of \eqnok{gen.60.a}, \eqnok{gen.60.b}, \eqnok{gen.60.c}, and
\eqnok{gen.12} can in turn be bounded by a finite quantity that
depends only on the problem data and not on $\mu$.

Returning to \eqnok{gen.10}, and using \eqnok{hatE}, we have that
\beqa
\nonumber
\bmat{c} w_{\hat{U}} \\ w_{\hat{V}} \\ y_U \emat &=&
-(Q + \hE_{11})^{-1} \hE_{12} y_V + 
(Q + \hE_{11})^{-1} \bmat{c} 
\hat{U}^T r_1 \\ \hat{V}^T r_1 \\ \Sigma \hat{U}^T r_3 + U^T r_4  \emat \\
\nonumber
&=& O(\| \hE_{12}\| \| y_V\|) + O(\| r_1 \| + \| r_3 \| + \| r_4 \|) \\
\label{gen.15}
&=& O(\mu) \| y_V \| +  O(\| r_1 \| + \| r_3 \| + \| r_4 \|).
\eeqa
Meanwhile, for the second block row of \eqnok{gen.10}, we obtain
\beq \label{gen.16}
y_V = -(\tilde{M}_{VV} + \hE_{22})^{-1} \hE_{21} \bmat{c}
w_{\hat{U}} \\ w_{\hat{V}} \\ y_U \emat +
(\tilde{M}_{VV} + \hE_{22})^{-1} V^T r_4.
\eeq
Since from \eqnok{mtilinv}, \eqnok{hatE}, and \eqnok{gen.0.1}, we have
\[
(\tilde{M}_{VV} + \hE_{22})^{-1}
= (I + \tilde{M}_{VV}^{-1}  \hE_{22})^{-1} \tilde{M}_{VV}^{-1}  
= (I + \dbu/\mu)  \tilde{M}_{VV}^{-1}  = O(\mu^{-1}),
\]
it follows  from \eqnok{gen.16} and \eqnok{hatE} that
\[
y_V = O(\mu^{-1}) O(\mu) \left\|
\bmat{c} w_{\hat{U}} \\ w_{\hat{V}} \\ y_U \emat \right\| 
+ O(\mu^{-1}) O(\| r_4 \|).
\]
By substituting from \eqnok{gen.15}, we obtain
\[
\| y_V \| = O(\mu) \| y_V\| +
O(\| r_1 \| + \| r_3 \| + \| r_4 \|) + O(\| r_4 \| / \mu),
\]
and therefore
\[
\| y_V \| = O(\| r_1 \| + \| r_3 \| + \| r_4 \| / \mu),
\]
as claimed. The estimate for $(w_{\hat{U}}, w_{\hat{V}}, y_U)$ is
obtained by substituting into \eqnok{gen.15}.
\end{proof}

The conditions \eqnok{gen.0} need a little explanation.  For the
typical value $\bu=10^{-16}$, the minimum value of the quantity
$\max(\mu^{1/3}, \bu/\mu)$ is $10^{-4}$, achieved at
$\mu^{-12}$. Moreover, we have $\max(\mu^{1/3}, \bu/\mu) \le 10^{-2}$
only for $\mu$ in the range $[10^{-14}, 10^{-6}]$.  It would seem,
then, that the problem would need to be quite well conditioned for
\eqnok{gen.0.3} and \eqnok{gen.0.2} to hold and that $\mu$ may have
to become quite small before the results apply. We note, however, that
the purpose of the bounds \eqnok{gen.0.3} and \eqnok{gen.0.2} is to
ensure that the inverse of $Q+\hE_{11}$ can be bounded independently of
$\mu$, and that for this purpose they are quite conservative. That is,
we would expect to find that $\| (Q+\hE_{11})^{-1}\| $ is not too much
larger than the norm of the inverse of the corresponding exact matrix
(the first term on the right-hand side of \eqnok{gen.10b}) for $\mu$
not much less than the smallest eigenvalues of $\Sigma$ and $\hat{V}^T
\cL_{zz}(z^*,\lambda^*) \hat{V}$.

The requirement that $\bu/\mu$ and $\mu$ both be small in
\eqnok{gen.0} may not seem to sit well with expressions such as
$O(\mu)$ and $O(\mu^2)$, which crop up repeatedly in the analysis and
which assert that certain bounds hold ``for all sufficiently small
$\mu$.''  
As noted in the preceding paragraph, this requirement implies that the
analysis holds for $\mu$ in a certain range, or ``window,'' of
values. Similar windows are used in the analysis of
S.~Wright~\cite{Wri93c,Wri94a_rev,Wri96b}, and
M.~H.~Wright~\cite{MWri98a}, and numerical experience indicates that
such a window does indeed exist in most practical cases. We expect the
same to be true of the problem and algorithms discussed in this paper.

At this point, we assemble the assumptions that are made in the
remainder of the paper into a single catch-all assumption.

\begin{assumption} \label{ass:base}
\mbox{}
\begin{itemize}
\item[(a)] $z^*$ is a solution of \eqnok{nlp}, so that the condition
\eqnok{kkt} holds. The MFCQ condition \eqnok{mfcq}, the second-order
conditions \eqnok{2os.L}, \eqnok{2os.1}, and the SC condition
\eqnok{sc} are satisfied at this solution.  The current iterate
$(z,\lambda,s)$ of the PDIP algorithm satisfies the conditions
\eqnok{pdip.central}, and the right-hand side modification $\tpert$
satisfies \eqnok{tpert.est}.
\item[(b)] The quantities $\mu$, $\bu$ \eqnok{def.bu},
$\cL_{zz}(z^*,\lambda^*)$, $\Sigma$, and $\hat{V}$ \eqnok{gen.3}
satisfy the conditions \eqnok{gen.0}.
\end{itemize}
\end{assumption}

From our observations following \eqnok{def.DLS}, we have under this
assumption that 
\beq  \label{dBN.est}
D_{\cB} = O(\mu), \sgap
D_{\cB}^{-1} = O(\mu^{-1}), \sgap
D_{\cN} = O(\mu^{-1}), \sgap
D_{\cN}^{-1} = O(\mu).
\eeq

Our next result considers a perturbed form of the system
\eqnok{gen:aug}, with a general right-hand side.  By eliminating one
component to obtain the form \eqnok{gen:augB}, we can apply
Theorem~\ref{th:gen1} to obtain estimates of the dependence of the
solution on the right-hand side components.
\begin{theorem} \label{th:gen2}
Suppose that Assumption~\ref{ass:base} holds.
Consider the linear system
%
%
\beqa
\nonumber
& \bmat{ccc} \cL_{zz}(z,\lambda) + E_{11} & \nablagB (z) + E_{12} & \nablagN (z) + E_{13} \\ 
\nablagB (z)^T + E_{21} & -D_{\cB} + E_{22} & E_{23} \\ 
\nablagN (z)^T + E_{31} & E_{32} & -D_{\cN} + E_{33}
\emat
\bmat{c}  w \\ y \\ q  \emat \\
\label{gen.30} 
& = 
\bmat{c} r_5 \\ \nablagB(z^*)^T r_6 + r_7 \\ r_8 \emat, 
\eeqa
where
\begin{subequations} \label{esizes}
\beqa 
\label{esizes.1}
& E_{11} = \dbu / \mu, \sgap E_{33} = \dbu / \mu^2, \\
\label{esizes.2}
& E_{12}, E_{21}, E_{22} = \dbu, \sgap E_{13}, E_{31}, E_{23}, E_{32} =
\dbu/\mu.
\eeqa
\end{subequations}
Then the step $(w,y,q)$ satisfies the following estimates:
\beqas
(U^Ty,w) &=& O(\|r_5\| + \|r_6 \| + \| r_7 \| + \mu \| r_8 \|), \\
V^Ty &=& O(\| r_5 \| + \| r_6 \| + \| r_7 \|/\mu + 
(\dbu  / \mu + O(\mu)) \| r_8 \|), \\
q &=& O(\mu) \left[ \| r_5 \| + \| r_6 \| + \| r_8 \| \right]
+ (\dbu/\mu + O(\mu)) \| r_7 \|.
\eeqas
\end{theorem}
\begin{proof}
%
  From \eqnok{dBN.est} and the assumed bound \eqnok{esizes.1} on the
  size of $E_{33}$, we have that
\beqa
\nonumber
\lefteqn{ (-D_{\cN} + E_{33})^{-1} } \\
\label{gen.34}
& = & -(I - D_{\cN}^{-1} E_{33} )^{-1} D_{\cN}^{-1} =
(I + O(\mu) \dbu/\mu^2 ) O(\mu) = O(\mu).
\eeqa
By eliminating $q$ from \eqnok{gen.30}, we obtain the reduced system
\[
\bmat{cc} H (z,\lambda) + \tE_{11} & \nablagB(z) + \tE_{12} \\
\nablagB(z)^T + \tE_{21} & -D_{\cB} + \tE_{22} \emat
\bmat{c} w \\ y \emat =
\bmat{c} r_5 + O(\mu ) \| r_8 \| \\ 
\nablagB(z^*)^T r_6 + r_7 + \dbu \| r_8 \| \emat,
\]
where from \eqnok{gen.0} and \eqnok{Hdef}, we obtain
\beqas
\tE_{11} &=& E_{11} - (\nablagN(z)+E_{13})(-D_{\cN} + E_{33})^{-1}
(\nablagN(z)^T +E_{31}) -\nablagN(z) D_{\cN}^{-1}\nablagN(z)^T \\
& = & \dbu/\mu + O(\mu), \\
\tE_{12} &=& E_{12} - (\nablagN(z)+E_{13})(-D_{\cN} + E_{33})^{-1}
E_{32} = \dbu + O(1) O(\mu) \dbu/\mu = \dbu, \\
\tE_{21} &=& E_{21} - E_{23} (-D_{\cN} + E_{33})^{-1}
(\nablagN(z)^T +E_{31}) = \dbu + (\dbu/\mu)  O(\mu) O(1) = \dbu, \\
\tE_{22} &=& E_{22} - E_{23}  (-D_{\cN} + E_{33})^{-1} E_{32} =
\dbu + (\dbu/\mu)^2 O(\mu) = \dbu.
\eeqas
These perturbation matrices satisfy the assumptions of
Theorem~\ref{th:gen1}, which can be applied to give
\begin{subequations} \label{gen.33}
\beqa \label{gen.33a}
(U^Ty , \hat{V}^Tw, \hat{U}^Tw) &=& O (\| r_5\| + \| r_6 \| + \| r_7 \| +
\mu \| r_8 \|), \\
V^Ty &=& O(\| r_5 \| + \|r_6 \| + \| r_7 \| / \mu ) + 
(\dbu/\mu + O(\mu)) \| r_8 \|.
\eeqa
\end{subequations}
From the last block row of \eqnok{gen.30}, and using \eqnok{gen.0},
\eqnok{gen.34}, \eqnok{gen.33}, we obtain
\beqas
q &=& (-D_{\cN} + E_{33})^{-1} \left[ r_8 - (\nablagN(z)^T + E_{31}) w -
E_{32} y \right] \\
&=& O(\mu) \left[ \| r_8 \| + \| w \| + (\dbu/\mu) \| y \| \right] \\
&=& O(\mu) \left[ \| r_5 \| + \| r_6 \| + \| r_7 \| + \| r_8 \| \right] + \\
&& \hspace*{.5in} \dbu \left[ \| r_5 \| + \| r_6 \| +
\| r_7 \| / \mu + ( \dbu / \mu + O(\mu) ) \| r_8 \| \right] \\
&=& O(\mu) \left[ \| r_5 \| + \| r_6 \| + \| r_8 \|  \right] +
(\dbu / \mu + O(\mu)) \| r_7 \|.
\eeqas
\end{proof}

An estimate for the solution of the exact system \eqnok{pdip.orig}
follows almost immediately from this result. This is the key technical
result used by Ralph and Wright~\cite{RalW96,RalW96b} to prove
superlinear convergence of PDIP algorithms for convex programming
problems. The result below, however, does not require a convexity
assumption.
\begin{corollary} \label{co:exact}
Suppose that Assumption~\ref{ass:base}(a) holds.
Then the (exact) solution
$(\Dz, \Dl, \Ds)$ of the system \eqnok{pdip.orig}  satisfies
\beq \label{est:exact}
(\Dz, \Dl, \Ds) = O(\mu).
\eeq
\end{corollary}
\begin{proof}
Note first that Assumption~\ref{ass:base}(b) holds trivially in
this case for $\mu$ sufficiently small, because our assumption of
exact computations is equivalent to setting $\bu=0$.
We prove the result by identifying the system \eqnok{gen:aug} with
\eqnok{gen.30} and then applying Theorem~\ref{th:gen2}.


For the right-hand side, we note first that, because of smoothness of
$g$, Taylor's theorem, the definition \eqnok{Bdef} of $\cB$, and
Theorem~\ref{th:est3},
\beqa
\nonumber
g_{\cB}(z) &=& g_{\cB}(z^*) + \nablagB(z^*)^T (z-z^*) + O(\| z-z^*\|^2) \\
\label{gBsize}
&=& \nablagB(z^*)^T (z-z^*) + O(\mu^2).
\eeqa
We now identify the right-hand side of \eqnok{gen:aug} with \eqnok{gen.30}
by setting
\beqas
r_5 &=& -\cL_z(z,\lambda), \\
r_6 &=&  (z-z^*), \\
r_7 &=& -\nablagB(z^*)^T (z-z^*) - \gB(z) + \LB^{-1} t_{\cB}, \\
r_8 &=& -\gN(z) + \LN^{-1} \tpert_{\cN}.
\eeqas
The sizes of these vectors can be estimated
by using \eqnok{pdip.central}, Lemma~\ref{lem:est2}, \eqnok{gBsize}, 
Theorem~\ref{th:est3}, and the assumption \eqnok{tpert.est} on
the size of $\tpert$  to obtain
\beq \label{rests.2}
r_5 = O(\mu), \sgap
r_6 = O(\mu), \sgap
r_7 = O(\mu^2), \sgap
r_8 = O(1).
\eeq
(By choosing $r_6$ and $r_7$ in this way, we ensure that the terms 
involving $\| r_7 \|/\mu$ in  the estimates of the solution components
in Theorem~\ref{th:gen2} are not grossly larger than the other
terms in these expressions.)
We complete the identification of \eqnok{gen:aug} with \eqnok{gen.30}
by setting all the perturbation matrices $E_{11}, E_{12}, \dots,
E_{33}$ to zero and by identifying the solution vector components
$\Dz$, $\Dl_{\cB}$, and $\Dl_{\cN}$ with $w$, $y$, and $q$,
respectively.  By directly applying Theorem~\ref{th:gen2},
substituting the estimates \eqnok{rests.2}, and setting $\dbu=0$ (since
we are assuming exact computations), we have that
\[
(U^T \Dl_{\cB}, \Dz) = O(\mu), \sgap
V^T \Dl_{\cB} = O(\mu), \sgap
\Dl_{\cN} = O(\mu).
\]

To show that the remaining solution component $\Ds$ of
\eqnok{pdip.orig} is also of size $O(\mu)$, we write the second block
row in \eqnok{pdip.orig} as
\[
\Ds = -(g(z) + s) - \nablag(z)^T \Dz,
\]
from which the desired estimate follows immediately by substituting
from \eqnok{pdip.central.2} and $\Dz = O(\mu)$.
\end{proof}

The next result uses Theorem~\ref{th:gen2} to compare perturbed and
exact solutions of the system of the system \eqnok{gen:aug}.
\begin{corollary} \label{co:errp}
Suppose that Assumption~\ref{ass:base} holds.
Let $(w,y,q)$ be obtained from the
following perturbed version of \eqnok{pdip.aug}:
\beqa
\nonumber & 
\bmat{ccc} 
\cL_{zz}(z,\lambda) + E_{11} & \nablagB (z) + E_{12} & \nablagN (z) + E_{13} \\ 
\nablagB (z)^T + E_{21} & -D_{\cB} + E_{22} & E_{23} \\ 
\nablagN (z)^T + E_{31} & E_{32} & -D_{\cN} + E_{33}
\emat \bmat{c} w \\ y \\ q \emat \\
\label{errp.1}
& \hspace*{1in} =
\bmat{c} -\cL_z(z,\lambda) + f_1 \\  
-\gB(z) + \LB^{-1} \tpert_{\cB} + f_2 \\
-\gN(z) + \LN^{-1} \tpert_{\cN} + f_3
\emat,
\eeqa
where $E_{ij}$, $i,j=1,2,3$, satisfy the conditions \eqnok{esizes} and
$f_1$, $f_2$, and $f_3$ are all of size $\dbu$.  Then if $(\Dz, \Dl,
\Ds)$ is the (exact) solution  of the system \eqnok{pdip.orig}, we have
\beqas
(\Dz-w, U^T(\DlB-y)) &=& \dbu, \\
V^T(\DlB - y) &=& \dbu/\mu, \\
\DlN -q &=& \dbu.
\eeqas
\end{corollary}
\begin{proof}
By combining \eqnok{errp.1} with \eqnok{gen:aug}, we obtain
\beqa
\nonumber & 
\bmat{ccc} 
\cL_{zz}(z,\lambda) + E_{11} & \nablagB (z) + E_{12} & \nablagN (z) + E_{13} \\ 
\nablagB (z)^T + E_{21} & -D_{\cB} + E_{22} & E_{23} \\ 
\nablagN (z)^T + E_{31} & E_{32} & -D_{\cN} + E_{33}
\emat \bmat{c} w -\Dz \\ y -\DlB \\ q - \DlN  \emat \\
\label{errp.2} &
\hspace*{1in} =
\bmat{c} f_1 \\ f_2 \\ f_3 \emat -
\bmat{ccc} E_{11} & E_{12} & E_{13} \\ E_{21} & E_{22} & E_{23} \\
E_{31} & E_{32} & E_{33} \emat 
\bmat{c} \Dz \\ \DlB \\ \DlN \emat.
\eeqa
From the bounds on the perturbations $E$ in \eqnok{esizes} and the
result of Corollary~\ref{co:exact}, we have for the
right-hand side of this expression that
\beqa
\nonumber
\bmat{c} r_5 \\ r_7 \\ r_8 \emat & \defeq  &
\bmat{c} f_1 \\ f_2 \\ f_3 \emat -
\bmat{ccc} E_{11} & E_{12} & E_{13} \\ E_{21} & E_{22} & E_{23} \\
E_{31} & E_{32} & E_{33} \emat 
\bmat{c} \Dz \\ \DlB \\ \DlN \emat \\
\label{bush.sux1}
&=& \bmat{c} 
\dbu + (\dbu/\mu) \mu + \dbu \mu + (\dbu/\mu) \mu \\
\dbu + \dbu \mu + \dbu \mu + (\dbu/\mu) \mu \\
\dbu + (\dbu/\mu) \mu + (\dbu/\mu) \mu + (\dbu/\mu^2) \mu
\emat = \bmat{c} \dbu  \\ \dbu \\ \dbu/\mu \emat.
\eeqa
Using these estimates, we can simply apply Theorem~\ref{th:gen2}  to
\eqnok{errp.2} (with $r_6=0$) to obtain the result.
\end{proof}

For later reference, we show how the estimates of
Corollary~\ref{co:errp} can be modified when the perturbations have a
special form. Suppose that
\beq \label{E23f2}
E_{23} = 0, \sgap E_{33} = \dbu/\mu, \sgap f_2 = U f_2^U + O(\mu^2), \;\; 
\mbox{where} \; f_2^U = \dbu,
\eeq
where $U$ is the matrix from \eqnok{gen.3}. Instead of setting
$r_6=0$ as in the proof above, we set
\[
r_6 = \hat{U} \Sigma f_2^U = \dbu
\]
(using \eqnok{gen.3} to obtain an $r_6$ for which $\nablagB(z^*)^T r_6
= U f_2^U$).  By modifying \eqnok{bush.sux1} to account for the
remaining perturbations, we can identify \eqnok{errp.2} with
\eqnok{gen.30} by setting
\beqa
\nonumber
\bmat{c} r_5 \\ r_7 \\ r_8 \emat & \defeq  &
\bmat{c} f_1 \\ f_2 - U f_2^U \\ f_3 \emat -
\bmat{ccc} E_{11} & E_{12} & E_{13} \\ E_{21} & E_{22} & E_{23} \\
E_{31} & E_{32} & E_{33} \emat 
\bmat{c} \Dz \\ \DlB \\ \DlN \emat \\
\label{bush.sux2}
&=& \bmat{c} 
\dbu + (\dbu/\mu) \mu + \dbu \mu + (\dbu/\mu) \mu \\
O(\mu^2) + \dbu \mu + \dbu \mu \\
\dbu + (\dbu/\mu) \mu + (\dbu/\mu) \mu + (\dbu/\mu^2) \mu
\emat = \bmat{c} \dbu  \\ O(\mu^2) \\ \dbu/\mu \emat.
\eeqa
Using these modified right-hand side estimates, we can apply
 Theorem~\ref{th:gen2} to obtain the following improved bound on 
one of the components:
\beq \label{improve.Vdl}
V^T(\DlB - y) = O(\mu).
\eeq
The bounds on the other components remain unchanged.

We emphasize that the conditions \eqnok{pdip.central}, and in
particular \eqnok{pdip.central.3}, are critical to the results of this
and all the remaining sections of the paper. These conditions enable
Lemma~\ref{lem:est2}, which in turn
enable us to assert that the diagonals of $D_{\cB}$ all have size
$\Theta(\mu)$ while those of $D_{\cN}$ all have size
$\Theta(\mu^{-1})$ (see \eqnok{dBN.est}). This neat classification of the diagonals of $D$
into two categories drives all the subsequent analysis. The motivation
for conditions like \eqnok{pdip.central} in the literature for
path-following methods (with exact steps) is not unrelated: It allows
us to obtain bounds on the steps and to show that we can move a
significant distance along this direction while ensuring that
\eqnok{pdip.central} continues to be satisfied at the new iterate.
(See, for example, \cite[Chapters~5 and 6]{IPPD96} and its bibliography
for the case of linear programming and \cite{SJW32,RalW96,RalW96b}
for the case of nonlinear convex programming.) In the analysis above,
we obtain bounds on the {\em errors} (rather than the steps
themselves) when perturbation terms of a certain structure appear in
the matrix and right-hand side.

Many practical implementations of path-following methods for linear
programming do not explicitly check that the condition
\eqnok{pdip.central.3} is satisfied by the calculated iterates (see,
for example, \cite{Meh92a} and \cite{PCx99}). However, the heuristics
for ``stepping back'' from the boundary of the nonnegative orthant by
a small but significant quantity are motivated by this condition, and
it is observed to hold in practice on all but the most recalcitrant
problems. 

\section{The Condensed System} \label{sec:pdip.cond}

Here we consider an algorithm in which the condensed linear system
\eqnok{pdip.cond} is formed and solved to obtain $\Dz$, and the
remaining step components $\Dl$ and $\Ds$ are recovered from
\eqnok{pdip.orig}. We obtain expressions for the errors in the
calculated step $(\wDz, \wDl, \wDs)$ and discuss the effects of these
errors on certain measures of step quality. We also derive conditions
under which the Cholesky factorization applied to \eqnok{pdip.cond} is
guaranteed to run to completion.

Formally, the complete procedure can be described as follows:

\medskip

\btab
\> {\bf procedure condensed} \\
\> {\bf given} the current iterate $(z,\lambda,s)$ \\
\>\> form the coefficient matrix and right-hand side for \eqnok{pdip.cond}; \\
\>\> solve \eqnok{pdip.cond} using a backward stable algorithm to obtain $\Dz$; \\
\>\> set $\Dl = D^{-1} [ g(z) - \Lambda^{-1} \tpert + \nablag(z)^T \Dz ]$; \\
\>\> set $\Ds = -(g(z) + s) - \nablag(z)^T \Dz$.
\etab

\medskip

We have used the  definition \eqnok{def.DLS} of the matrix $D$.
For convenience, we restate the system \eqnok{pdip.cond} here as follows:
\beq \label{pdip.cond2}
\left[ \cL_{zz} (z,\lambda) + \nablag(z) D^{-1} \nablag(z)^T 
\right] \Dz = -\cL_z (z,\lambda) - \nablag(z) D^{-1}
[g(z) - \Lambda^{-1} \tpert ].
\eeq

\medskip

Note that this procedure requires evaluation of $D^{-1} = S^{-1}
\Lambda$, rather than $D$ itself.

\subsection{Quantifying the Errors} \label{sec:cond:quant}

When implemented in finite-precision arithmetic, 
solution of \eqnok{pdip.cond2}
gives rise to errors of three types:
\bi
\item[-] cancellation in evaluation of the matrix and right-hand side;
\item[-] roundoff errors in evaluation of the matrix and right-hand side;
\item[-] roundoff errors that accumulate during the process of
factoring the matrix and using triangular substitutions to obtain the
solution. 
\ei

Cancellation may be an issue in the evaluation of the nonlinear
functions $\cL_{zz}(z,\lambda)$, $\cL_z(z,\lambda)$, $g(z)$, and
$\nablag(z)$, because intermediate terms computed during the additive
evaluation of these quantities may exceed the size of the final result
(see Golub and Van Loan~\cite[p.~61]{GolV96}).
The intermediate terms generally contain rounding error (which occurs whenever
real numbers are represented in finite precision).
Cancellation becomes a
significant phenomenon whenever we take a difference of two nearly
equal quantities, since the error in the computed result due to
roundoff in the two arguments may be large relative to the size of the
result.  If, as we can reasonably assume, intermediate quantities in
the calculations of our right-hand sides remain bounded, the absolute
size of the errors in the result is $\dbu$.
In the case of $\cL_z(z,\lambda)$ and $\gB(z)$, the final result in
exact arithmetic has size $O(\mu)$, so that the error $\dbu$ takes on
a large relative significance for small values of $\mu$. This fact
causes the error bound in some components of the solution to be larger
than in others, as we see in \eqnok{cond.p4c} below.  In summary, the
computed versions of the quantities discussed above differ from their
exact values in the following way:
\begin{subequations} \label{cancel}
\beqa
\label{cancel.1}
\makebox{\rm computed } \cL_{zz}(z,\lambda)  \leftarrow 
\cL_{zz}(z,\lambda) + \bar{F}, \\
\label{cancel.2}
\makebox{\rm computed } \cL_z (z,\lambda)   \leftarrow 
\cL_z (z,\lambda) + \bar{f}, \\
\label{cancel.3}
\makebox{\rm computed } \nablag(z)  \leftarrow 
\nablag(z) + F = \bmat{c} \nablagB(z) \\ \nablagN(z) \emat + 
\bmat{c} F_{\cB} \\ F_{\cN} \emat, \\
\label{cancel.4}
\makebox{\rm computed } g(z) \leftarrow 
g(z) + f = \bmat{c} \gB(z) \\ \gN(z) \emat +
\bmat{c} f_{\cB} \\ f_{\cN} \emat,
\eeqa
\end{subequations}
where $\bar{F}$, $\bar{f}$, $F$, and $f$ are all of size $\dbu$.
Earlier discussion of cancellation in similar contexts can be found in
the papers of S.~Wright~\cite{Wri93c,Wri94a_rev,Wri96b} and
M.~H.~Wright~\cite{MWri98a}.

The second source of error is evaluation of the matrix $D^{-1}$. From
the model \eqnok{def.bu} of floating-point arithmetic and the
estimates \eqnok{eq:14} of Lemma~\ref{lem:est2}, we have that
\begin{subequations} \label{cancelD}
\beqa \label{cancelD.B}
\makebox{\rm computed } D_{\cB}^{-1}  & \leftarrow & (D_{\cB} + G_{\cB})^{-1},
 \sgap G_{\cB} = \mu \dbu, \\
\label{cancelD.N}
\makebox{\rm computed } D_{\cN}^{-1} & \leftarrow & (D_{\cN} + G_{\cN})^{-1}, 
\sgap G_{\cN} =  \dbu / \mu,
\eeqa
\end{subequations}
where $G_{\cB}$ and $G_{\cN}$ are both diagonal matrices that can be
composed into a single diagonal matrix $G$. 

Third, we account for the error in forming the matrix and right-hand
side of \eqnok{pdip.cond2} from the computed quantities described in
the last two paragraphs. Since we are now dealing with floating-point
numbers, the model \eqnok{def.bu} applies; that is, any additional
errors that arise during the combination of these floating-point
quantities have size $\bu$ relative to the size of the result of the
calculation.  Since the norm of the coefficient matrix is of size
$O(\mu^{-1})$ while the right-hand side has size $O(1)$ (see
\eqnok{pdip.central}), we represent these errors by a matrix $\hat{F}$
of size $\dbu/\mu$ and a vector $\hat{f}$ of size $\dbu$.

Finally, we account for the error that arises in the application of a
backward-stable method to solve \eqnok{pdip.cond2}. Specifically, we
assume that the method yields a computed solution that is the exact
solution of a nearby problem whose data contains relative
perturbations of size $\bu$.  The absolute sizes of these terms would
therefore be $\dbu/\mu$ in the case of the matrix and $\dbu$ in the
case of the right-hand side. Since these errors are the same size as
those discussed in the preceding paragraph, we incorporate them into
the matrix $\hat{F}$ and the vector $\hat{f}$.

Summarizing, we find that the computed solution $\wDz$ of
\eqnok{pdip.cond2} satisfies the following system:
\beqa 
\label{pdip.cond2.err}
\lefteqn{
\left[ \cL_{zz} (z,\lambda) + \bar{F} + (\nablag(z)+F) (D+G)^{-1} (\nablag(z)+F)^T  
+ \hat{F} \right] \wDz } \\
\nonumber
& = & -\cL_z (z,\lambda)  - \bar{f} - (\nablag(z)+F) (D+G)^{-1}
[g(z) +f - \Lambda^{-1} \tpert ] + \hat{f},
\eeqa
where the perturbation terms $\bar{F}$, $F$, $\hat{F}$, $G$,
$\bar{f}$, $\hat{f}$, and $f$ are described in the paragraphs above.
By ``unfolding'' this system and using the partitioning of $F$, $G$,
and $f$ defined in \eqnok{cancel} and \eqnok{cancelD}, we find that
$\wDz$ also satisfies the following system, for some vectors $y$ and
$q$:
\beqa \label{cond.p3} 
& \bmat{ccc} \cL_{zz}(z,\lambda) +  \bar{F} + \hat{F} 
& \nablagB (z) + F_{\cB} 
& \nablagN (z) + F_{\cN} \\ 
\nablagB (z)^T + F_{\cB}^T & 
-D_{\cB} - G_{\cB} & 0 \\
\nablagN (z)^T + F_{\cN}^T & 0 & -D_{\cN} - G_{\cN}
\emat \bmat{c} \wDz \\ y \\ q \emat \\
\nonumber
& \hspace*{1in} = \bmat{c} -\cL_z(z,\lambda) - \bar{f} + \hat{f} \\ 
-\gB(z) + \LB^{-1} \tpert_{\cB} - f_{\cB} \\
-\gN(z) + \LN^{-1} \tpert_{\cN} - f_{\cN} \emat.
\eeqa
This system has precisely the form of \eqnok{errp.1} (in particular,
the perturbation matrices satisfy the appropriate bounds). Hence, by a
direct application of Corollary~\ref{co:errp}, we conclude that
\begin{subequations} \label{cond.p4}
\beqa \label{cond.p4a}
\Dz - \wDz & = & \dbu, \\
 \label{cond.p4b}
U^T( \Dl_{\cB} - y) &= & \dbu, \\
\label{cond.p4c}
V^T( \Dl_{\cB} - y) &= & \dbu / \mu.
\eeqa
\end{subequations}

We return now to the recovery of the remaining solution components
$\wDl$ and $\wDs$ in the procedure {\bf condensed}. 
We have from Assumption~\ref{ass:base} together with
\eqnok{pdip.central.2},
Lemma~\ref{lem:est2},
\eqnok{est:exact},
\eqnok{cond.p4a},
\eqnok{cancelD.B},
\eqnok{gen.0}, and \eqnok{dBN.est} that
%
\begin{subequations} \label{cond.p8} 
\beqa \label{cond.p8a}
& \gB(z) = r_g(z,s)_{\cB} - s_{\cB} = O(\mu), \sgap 
\LB^{-1} = \Theta(1), \sgap
\wDz = \Dz + \dbu = O(\mu), \\
\label{cond.p8b}
& (D_{\cB} + G_{\cB})^{-1} 
= (I + D_{\cB}^{-1} G_{\cB}) D_{\cB}^{-1} = (I + \dbu)^{-1} O(\mu^{-1}) =
O(\mu^{-1}).
\eeqa
\end{subequations}
Since $\tpert=O(\mu^2)$, we have from our model \eqnok{def.bu} that
the floating-point version of the calculation of $\wDl_{\cB}$ in the
procedure {\bf condensed} satisfies the following:
\[
\wDl_{\cB} = (D_{\cB} + G_{\cB})^{-1} \left[
\gB(z) + f_{\cB} - \LB^{-1} \tpert_{\cB} +
(\nablagB(z) + F_{\cB})^T \wDz + \mu \dbu \right] + \dbu.
\]
(The final term $\dbu$ arises from \eqnok{def.bu} because our best
estimate of the quantity in the brackets at this point of the analysis
is $O(\mu)$, so from \eqnok{cond.p8b} the result has size $O(1)$.)
Meanwhile, we have from the second block row of \eqnok{cond.p3} that 
\[
y = (D_{\cB} + G_{\cB})^{-1} \left[ 
\gB(z) + f_{\cB} - \LB^{-1} \tpert_{\cB} +
(\nablagB(z) + F_{\cB})^T \wDz \right].
\]
By a direct comparison of these two expressions, and using $(D_{\cB} +
G_{\cB})^{-1} = O(\mu)$, we find that
\beq \label{cond.p9}
\wDl_{\cB} - y = \dbu.
\eeq
By combining \eqnok{cond.p9} with \eqnok{cond.p4b} and
\eqnok{cond.p4c}, we find that
\beq \label{cond.p10}
U^T(\Dl_{\cB} - \wDl_{\cB}) = \dbu, \sgap
V^T(\Dl_{\cB} - \wDl_{\cB}) = \dbu / \mu.
\eeq

For the ``nonbasic'' part $\wDl_{\cN}$, we have from 
\eqnok{pdip.central.2},
Lemma~\ref{lem:est2},
\eqnok{est:exact},
\eqnok{cond.p4a},
\eqnok{cancelD.N},
\eqnok{gen.0}, and \eqnok{dBN.est} that
%
\begin{subequations} \label{cond.p11}
\beqa \label{cond.p11a}
& \gN(z) = O(1), \sgap \LN^{-1} = O(\mu^{-1}), \sgap \wDz = O(\mu), \\
 \label{cond.p11b}
& (D_{\cN} + G_{\cN})^{-1} = (I + D_{\cN}^{-1} G_{\cN})^{-1} D_{\cN}^{-1} =
 D_{\cN}^{-1} + \mu \dbu = O(\mu).
\eeqa
\end{subequations}
By using $\tpert_{\cN} = O(\mu^2)$ and applying the model
\eqnok{def.bu} to the appropriate step in the procedure {\bf
condensed}, we obtain
\[
\wDl_{\cN} = (D_{\cN} + G_{\cN})^{-1} \left[
\gN(z) + f_{\cN} - \LN^{-1} \tpert_{\cN} +
(\nablagN(z) + F_{\cN})^T \wDz + \dbu \right] + \mu \dbu.
\]
By comparing this expression with the corresponding exact formula,
which is
\[
\Dl_{\cN} = D_{\cN}^{-1} \left[
\gN(z) - \LN^{-1} \tpert_{\cN} + \nablagN(z)^T \Dz \right],
\]
and by using the bounds  \eqnok{cond.p11} and the fact that 
$f_{\cN}$ and $F_{\cN}$ have size $\dbu$, we obtain
\beqa
\nonumber
\wDl_{\cN} - \Dl_{\cN} &=&
\mu \dbu + (D_{\cN} + G_{\cN})^{-1} [ f_{\cN} + \dbu ] + \\
\nonumber
&& \sgap [ (D_{\cN} + G_{\cN})^{-1} - D_{\cN}^{-1} ]
[\gN(z) - \LN^{-1} \tpert_{\cN}] + \\
\nonumber
&& \sgap (D_{\cN} + G_{\cN})^{-1} (\nablagN(z) + F_{\cN})^T \wDz -
D_{\cN}^{-1} \nablagN(z)^T \Dz \\
\nonumber 
&=& \mu \dbu + (D_{\cN} + G_{\cN})^{-1} 
[ \nablagN(z)^T (\wDz - \Dz) + F_{\cN}^T \wDz ] \\
\nonumber 
&& \sgap [ (D_{\cN} + G_{\cN})^{-1} - D_{\cN}^{-1} ] \nablagN(z)^T \Dz \\
\label{cond.p13} 
&=& \mu \dbu.
\eeqa

Finally, for the recovered step $\wDs$, we have from
the last step of procedure {\bf condensed}, together with
\eqnok{pdip.central.2}, \eqnok{cancel.4}, \eqnok{cond.p8b},
and \eqnok{def.bu} that
\[
\wDs = -(g(z) + f + s) - (\nablag(z) + F)^T \wDz + \dbu,
\]
where the final term accounts for the rounding error \eqnok{def.bu}
that arises from accumulating the terms in the sum, which are all bounded. 
By substituting the expression for the exact
$\Ds$ together with the estimates \eqnok{cancel.4} and
\eqnok{cond.p8b} on the sizes of the perturbation terms, we obtain
\beqa
\nonumber
\wDs &=& -(g(z) + s) - \nablag(z)^T \Dz -f - \nablag(z)^T (\wDz - \Dz) -
F^T \wDz + \mu \dbu  \\
\label{cond.p14}
&=& \Ds + \dbu.
\eeqa

We summarize the results obtained so far in the following theorem.
\begin{theorem} \label{th:cond1}
Suppose that Assumption~\ref{ass:base} holds.
Then when
the step \\
$(\wDz, \wDl, \wDs)$ is calculated in a finite-precision
environment by using the procedure {\bf condensed} (and
where, in particular, a backward stable method is used to solve the
linear system for the $\wDz$ component), we have that
\begin{subequations} \label{cond.p16}
\beqa
\label{cond.p16a}
(\Dz - \wDz, U^T (\Dl_{\cB} - \wDl_{\cB}), \Ds - \wDs) &=& \dbu, \\
\label{cond.p16b}
V^T (\Dl_{\cB} - \wDl_{\cB}) &=&  \dbu / \mu, \\
\label{cond.p16c}
\Dl_{\cN} - \wDl_{\cN} &=& \mu \dbu.
\eeqa
\end{subequations}
\end{theorem}

This theorem extends the result of M.~H.~Wright~\cite{MWri98a} for
accuracy of the computed solution of the condensed system by relaxing
the LICQ assumption to MFCQ. When LICQ holds, the matrix $V$ is
vacuous, so the absolute error in all components is of size at most
$\dbu$. The higher accuracy \eqnok{cond.p16c} of the components
$\wDl_{\cN}$ (also noted in \cite{MWri98a}) does not contribute
significantly to the progress that can be made along the inexact
direction $(\wDz, \wDl, \wDs)$, in the sense of
Section~\ref{sec:cond:local}.

We return briefly to the case discussed immediately after
Corollary~\ref{co:errp}, in which the perturbations have the special
form \eqnok{E23f2}, using these results to show that the bound
\eqnok{cond.p16b} can be strengthened when $f_{\cB}$ satisfies
\beq \label{improved.vtf}
V^T f_{\cB} = O(\mu^2). 
\eeq
This case is of interest when the cancellation errors in computing
$\gB(z)$ are smaller than the estimate we made following
\eqnok{cancel.4}, possibly because of use of higher-precision
arithmetic or the fact that the computation did not require
differencing of quantities whose size is large relative to the final
result.  When \eqnok{improved.vtf} holds, we see by comparing
\eqnok{errp.1} with \eqnok{cond.p3} that
\[
E_{23} = 0, \sgap E_{33} = G_{\cN} = \dbu/\mu, \sgap
f_2 = U U^T f_{\cB} + O(\mu^2), \;\; \mbox{where $f_{\cB} = \dbu$}.
\]
Therefore, we deduce from \eqnok{improve.Vdl} that \eqnok{cond.p4c} can be
replaced by 
\[
V^T( \Dl_{\cB} - y) = O(\mu).
\]
Using \eqnok{cond.p9} and $\mu \gg \dbu$, we can therefore replace
\eqnok{cond.p16b} in this case by
\beq \label{cond.p16b.improved}
V^T (\Dl_{\cB} - \wDl_{\cB}) =  O(\mu). 
\eeq

\subsection{Termination of the  Cholesky Algorithm} 
\label{sec:cond:chol}

In deriving the estimate \eqnok{cond.p4}, we have assumed that a
backward stable algorithm is used to solve \eqnok{pdip.cond2}.
Because of \eqnok{2os.L}, \eqnok{2os.1}, and the SC condition, and the
estimates of the sizes of the diagonals of $D$ (from \eqnok{def.DLS}
and Lemma~\ref{lem:est2}), it is easy to show that the matrix in
\eqnok{pdip.cond2} is positive definite for all sufficiently small
$\mu$. The Cholesky algorithm is therefore an obvious candidate for
solving this system. However, the condition number of the matrix in
\eqnok{pdip.cond2} usually approaches $\infty$ as $\mu \downarrow 0$,
raising the possibility that the Cholesky algorithm may break down
when $\mu$ is small. A simple argument, which we now sketch, suffices
to show that successful completion of the Cholesky algorithm can be
expected under the assumptions we have used in our analysis so far.

We state first the following technical result. Since it is similar to
one proved by 
Debreu~\cite[Theorem~3]{Deb52}, its proof is omitted.
\begin{lemma} \label{lem:tech1}
Suppose that $M$ and $A$ are two matrices with the properties that $M$
is symmetric and
\[
A^Tx=0 \;\; \Rightarrow \;\; x^TMx \ge \alpha \| M \|  \|x \|^2,
\]
for some constant $\alpha>0$. Then for all $\mu$ such that 
\[
0 < \mu < \bar{\mu} \defeq \min \left( \frac{\alpha \|A \|^2}{4 \|M \|}, 
\frac{\|A \|}{\alpha \|M \|} \right),
\]
we have that 
\[
x^T (M + \mu^{-1} AA^T) x \ge \frac{\alpha}{2} \|x\|^2, \sgap
\makebox{\rm for all $x$.}
\]
\end{lemma} 

We apply this result to \eqnok{pdip.cond2} by setting
\beqas
M &=& \cL_{zz} (z,\lambda) + \nablagN(z) D_{\cN}^{-1} \nablagN(z)^T =
 \cL_{zz} (z,\lambda) +  O(\mu), \\
A &=& \mu^{1/2} \nablagB(z) D_{\cB}^{-1/2}
\eeqas
(where again we use \eqnok{def.DLS} and Lemma~\ref{lem:est2} to derive
the order estimates).  The conditions \eqnok{2os.L}, \eqnok{2os.1},
and strict complementarity ensure that this choice of $M$ and $A$
satisfies the assumptions of Lemma~\ref{lem:tech1}. 
The result then implies that the smallest singular value of the
matrix in \eqnok{pdip.cond2} is positive and of size
$\Theta(1)$ for all values of $\mu$ below a threshold that is also of
size $\Theta(1)$. Since $D = O(\mu^{-1})$, the largest eigenvalue of
this matrix is of size $O(\mu^{-1})$,  so we have
\beq \label{cond.p5}
{\rm cond}(\cL_{zz} (z,\lambda) + \nablag(z) D^{-1} \nablag(z)^T ) = O(\mu^{-1}).
\eeq
(An estimate similar to this is derived by
M.~H.~Wright~\cite[Theorem~3.2]{MWri98a}, under the LICQ assumption.)
It is known from a result of Wilkinson (cited by  Golub and
Van Loan~\cite[p.~147]{GolV96}) that the Cholesky algorithm runs to
completion if $q_n \dbu {\rm cond}(\cdot) \le 1$, where $q_n$ is a
modest quantity that depends polynomially on the dimension $n$ of the
matrix. By combining this result with \eqnok{cond.p5}, we conclude
that for the matrix in \eqnok{pdip.cond2}, we can expect completion of
the Cholesky algorithm whenever $\mu \gg \dbu$. That is, no new
assumptions need to be added to those made in deriving the results of
earlier sections.

We note that this situation differs a little from the case of linear
programming where, because second-order conditions are not applicable,
it is usually necessary to modify the Cholesky procedure to ensure
that it runs to completion (see \cite{Wri96b}).


\subsection{Local Convergence with Computed Steps}
\label{sec:cond:local}

We begin this section by showing how the quantities $r_f$, $r_g$, and
$\mu$ change along the computed step $(\wDz, \wDl, \wDs)$ obtained
from the finite-precision implementation of the procedure {\bf
  condensed}.  We compare these with the changes that can be expected
along the exact direction $(\Dz, \Dl, \Ds)$.  We then consider the
effects of these perturbations on an algorithm of the type in which
the iterates are expected to satisfy the conditions
\eqnok{pdip.central}. Rapidly convergent variants of these algorithms
for linear programming problems usually allow the values of $C$ and
$\gamma$ in these conditions to be relaxed, so that a near-unit step
can be taken. We address the following question: If similar
relaxations are allowed in an algorithm for nonlinear programming, are
near-unit steps still possible when the steps contain perturbations of
the type considered above?

We show in particular that for the computed search direction, the
maximum step length that can be taken without violating the
nonnegativity conditions on $\lambda$ and $s$ satisfies
\beq \label{cond.p27}
1-\walfmax = \dbu /\mu + O(\mu),
\eeq
while the reductions in pairwise products, $\mu$, $r_f$, and $r_g$ satisfy
\begin{subequations} \label{cond.p28}
\beqa
\label{cond.p28p}
(\lambda_i +  \alpha \wDl_i)(s_i+\alpha \wDs_i) &=&
(1-\alpha) \lambda_i s_i + \dbu + O(\mu^2), \;\; i=1,2,\dots,m, \\
\label{cond.p28a}
\mu(\lambda+ \alpha \wDl, s + \alpha \wDs) &=& 
(1-\alpha) \mu + \dbu + O(\mu^2), \\
\label{cond.p28b}
r_f(z+\alpha \wDz, \lambda + \alpha \wDl) &=& 
(1-\alpha) r_f + \dbu + O(\mu^2), \\
\label{cond.p28c}
r_g(z+\alpha \wDz, s + \alpha \wDs) &=& 
(1-\alpha) r_g + \dbu + O(\mu^2).
\eeqa
\end{subequations}
The corresponding maximum steplength for the {\em exact} direction
satisfies
\beq \label{cond.p25}
1-\alfmax = O(\mu),
\eeq
while the reductions in $r_f$, $r_g$, and $\mu$ satisfy
\begin{subequations} \label{cond.p26}
\beqa
\label{cond.p26p}
(\lambda_i +  \alpha \Dl_i)(s_i+\alpha \Ds_i) &=&
(1-\alpha) \lambda_i s_i + O(\mu^2), \;\; i=1,2,\dots,m, \\
\label{cond.p26a}
\mu(\lambda+ \alpha \Dl, s + \alpha \Ds) &=& 
(1-\alpha) \mu + O(\mu^2), \\
\label{cond.p26b}
r_f(z+\alpha \Dz, \lambda + \alpha \Dl) &=& 
(1-\alpha) r_f + O(\mu^2), \\
\label{cond.p26c}
r_g(z+\alpha \Dz, s + \alpha \Ds) &=& 
(1-\alpha) r_g + O(\mu^2).
\eeqa
\end{subequations}
Our proof of the estimates \eqnok{cond.p27} and \eqnok{cond.p28} is
tedious but not completely straightforward, and we have included it in
the Appendix.


It is clear from \eqnok{cond.p27} and \eqnok{cond.p28} that the
direction $(\wDz, \wDl, \wDs)$ makes good progress toward the solution
set $\cS$. If the actual steplength $\alpha$ is close to its maximum
value $\walfmax$, in the sense that
\beq \label{alf.4}
\walfmax - \alpha =  \dbu/\mu + O(\mu),
\eeq
we have by direct substitution in \eqnok{cond.p27} and
\eqnok{cond.p28} that
\beqas
\mu(\lambda+ \alpha \wDl, s + \alpha \wDs) &=& \dbu + O(\mu^2), \\
r_f(z+\alpha \wDz, \lambda + \alpha \wDl) &=&  \dbu + O(\mu^2), \\
r_g(z+\alpha \wDz, s + \alpha \wDs) &=& \dbu + O(\mu^2).
\eeqas
These formulae suggest that finite precision does not have an
observable effect on the quadratic convergence rate of the
underlying algorithm until $\mu$ drops below
about $\sqrt{\bu}$.
%
Stopping criteria for interior-point methods usually include a
condition such as $\mu \le 10^4 \bu$ or $\mu \le \sqrt{\bu}$ (see, for
example, \cite{PCx99}), so that $\mu$ is not allowed to become so small
that the assumption $\mu \gg \bu$ made in \eqnok{gen.0} is violated.

In making this back-of-the-envelope assessment, however, we have not
taken into account the approximate centrality conditions
\eqnok{pdip.central}, which must continue to hold (possibly in a
relaxed form) at the new iterate. These conditions play a central role
both in the analysis above and in the convergence analysis of the
underlying ``exact'' algorithms, and also appear to be important in
practice. Typically (see, for example, Ralph and
Wright~\cite{RalW96}), the conditions \eqnok{pdip.central} are relaxed
by allowing a modest increase in $C$ and a modest decrease in $\gamma$
on the rapidly convergent steps. We show in the next result that
enforcement of these relaxed conditions is not inconsistent with
taking a step length $\alpha$ that is close to $\walfmax$, so that
rapid convergence can still be observed even in the presence of finite-precision effects.
\begin{theorem} \label{th:cond2}
Suppose that Assumption~\ref{ass:base} holds.
Then when the step \\
$(\wDz, \wDl, \wDs)$ is calculated in a
finite-precision environment by using the procedure {\bf condensed}, 
there is a constant $\hat{C}$ such that for all $\tau \in
[0,1/2]$ and all $\alpha$ satisfying
\beq \label{cent.1}
\alpha \in [0,1-\hat{C} \tau^{-1} (\bu/\mu + \mu) ],
\eeq
the following relaxed form of the approximate centrality conditions
holds:
\begin{subequations} \label{pdip.centralr}
\beqa
\label{pdip.central.1r}
r_f(z+\alpha \wDz,\lambda+\alpha \wDl) & \le &
C (1+\tau) \mu(\lambda + \alpha \wDl, s + \alpha \wDs), \\
\label{pdip.central.2r}
r_g (z+\alpha \wDz,s + \alpha \wDs) & \le & C (1+\tau) \mu(\lambda + \alpha \wDl, s + \alpha \wDs), \\
\label{pdip.central.3r}
(\lambda_i + \alpha \wDl_i) (s_i + \alpha \wDs_i) & \ge & \gamma (1-\tau) 
\mu(\lambda + \alpha \wDl, s + \alpha \wDs), \\
\nonumber && \gap \makebox{\rm for all $i=1,2,\dots,m$,}
\eeqa
where $C$ is the constant from conditions \eqnok{pdip.central}.
\end{subequations}
Moreover, when we set $\alpha$ to its upper bound in \eqnok{cent.1},
we find that
\beq \label{cent.2}
\delta(z+\alpha \wDz, \lambda + \alpha \wDl) \le \tau^{-1} (\dbu + O(\mu^2)).
\eeq
\end{theorem}
\begin{proof}
From \eqnok{pdip.central} and \eqnok{cond.p28}, we have that
\beqas
\lefteqn{
\| r_f(z+\alpha \wDz,  \lambda + \alpha \wDl) \|} \\
&=& (1-\alpha) \| r_f \| + \dbu + O(\mu^2) \\
&\le & C (1-\alpha) \mu + \dbu + O(\mu^2) \\
&=&  C (1+\tau) (1-\alpha) \mu
- C\tau (1-\alpha) \mu  + \dbu + O(\mu^2) \\
&=&  C (1+\tau) \mu(\lambda + \alpha \wDl, s + \alpha \wDs) 
- C\tau (1-\alpha) \mu  + \dbu + O(\mu^2).
\eeqas
We deduce that the required condition \eqnok{pdip.central.1r} will
hold provided that
\[
\dbu + O(\mu^2) \le  C\tau (1-\alpha) \mu.
\]
Since by definition we have that $\dbu +O(\mu^2) \le \bar{C} (\bu +
\mu^2)$ for some positive constant $\bar{C}$, we find that a
sufficient condition for the required inequality is that
\[
(1-\alpha) \ge (\bar{C}/C) \tau^{-1}  (\bu / \mu + \mu),
\]
which is equivalent to \eqnok{cent.1} for an obvious definition of
$\hat{C}$. Identical logic can be applied to $\| r_g \|$ to yield a
similar condition on $\alpha$.

For the condition \eqnok{pdip.central.3r}, we have from
\eqnok{pdip.central} and \eqnok{cond.p28} that
\beqas
\lefteqn{ (\lambda_i + \alpha \wDl_i) (s_i + \alpha \wDs_i) } \\
&=& (1-\alpha) \lambda_i s_i + \dbu +  O(\mu^2) \\
& \ge & (1-\alpha) \gamma \mu + \dbu + O(\mu^2) \\
&=& \gamma (1-\tau) (1-\alpha) \mu + \gamma \tau (1-\alpha) \mu + \dbu + O(\mu^2) \\
& = &
\gamma (1-\tau) \mu(\lambda + \alpha \wDl, s + \alpha \wDs) + 
\gamma \tau (1-\alpha) \mu + \dbu + O(\mu^2).
\eeqas
Hence, the condition \eqnok{pdip.central.3r} holds provided that
\[
\gamma \tau (1-\alpha) \mu + \dbu + O(\mu^2) \ge 0.
\]
Similar logic can be applied to this inequality to derive a bound of
the type \eqnok{cent.1}, after a possible adjustment of $\hat{C}$.

Finally, we obtain \eqnok{cent.2} by substituting $\alpha = 1-\hat{C}
\tau^{-1} (\bu / \mu  + \mu)$ into \eqnok{cond.p28} and applying
Theorem~\ref{th:est3}. (Note that, despite the relaxation of the
centrality conditions \eqnok{pdip.centralr}, the result of
Theorem~\ref{th:est3} still holds; we simply modify the proof to
replace $C$ by $(3/2)C$ in \eqnok{pdip.central.1} and
\eqnok{pdip.central.2}, and $\gamma$ by $\gamma/2$ in
\eqnok{pdip.central.3}.)
\end{proof}

\section{The Augmented System} \label{sec:pdip.aug}

In this section, we consider the case in which the augmented system
\eqnok{pdip.aug} (equivalently, \eqnok{gen:aug}) is solved to obtain
$(\Dz, \Dl)$, while the remaining step component $\Ds$ is recovered
from \eqnok{pdip.orig}. The formal specification for this procedure is
as follows:

\medskip

\btab
\> {\bf procedure augmented} \\
\> {\bf given} the current iterate $(z,\lambda,s)$ \\
\>\> form the coefficient matrix and right-hand side for \eqnok{gen:aug}; \\
\>\> solve \eqnok{gen:aug} to obtain $(\Dz,\Dl)$; \\
\>\> set $\Ds = -(g(z) + s) - \nablag(z)^T \Dz$.
\etab

\medskip

Much of our work in analyzing the augmented system form
\eqnok{gen:aug} has already been performed in Section~\ref{accuracy};
the main error result is simply Corollary~\ref{co:errp}.  However, we
can apply this result only if the floating-point errors made in
evaluating and solving this system satisfy the assumptions of this
corollary. In particular, we need to show that the perturbation
matrices $E_{ij}$, $i,j=1,2,3$ in \eqnok{errp.1} satisfy the estimates
\eqnok{esizes}.

This task is not completely straightforward. Unlike the condensed and
full-system cases, it is not simply a matter of assuming that a
backward-stable algorithm has been used to solve the system
\eqnok{gen:aug}. The reason is that the largest terms in the
coefficient matrix in \eqnok{pdip.aug}---the diagonal elements in the
matrix $D_{\cN}$---have size $O(\mu^{-1})$. The usual analysis of
backward-stable algorithms represents the floating-point errors as a
perturbation of the entire coefficient matrix whose size is bounded by
$\dbu$ times the matrix norm---in this case, $\dbu/\mu$.  However,
Corollary~\ref{co:errp} requires some elements of the perturbation
matrix to be smaller than this estimate; in particular, the
submatrices $E_{12}$, $E_{21}$, and $E_{22}$ need to be of size
$\dbu$. Therefore, we need to look closely at the particular
algorithms used to solve \eqnok{gen:aug} to see whether they satisfy the
following condition.

\begin{condition} \label{cond.1}
The solution obtained by applying the algorithm in
 question to the system \eqnok{gen:aug} in floating-point arithmetic
 is the exact solution of a perturbed system in which the
 perturbations of the coefficient matrix satisfy the estimates
 \eqnok{esizes}, while the right-hand side is unperturbed.
\end{condition}

We focus on diagonal pivoting methods, which take a symmetric matrix
$T$ and produce a factorization of the form
\beq \label{aug:dpiv}
PTP^T = LYL^T,
\eeq
where $P$ is a permutation matrix, $L$ is unit lower triangular, and
$Y$ is block diagonal, with a combination of $1 \times 1$ and
symmetric $2 \times 2$ blocks. The best-known methods of this class
are due to Bunch and Parlett~\cite{BunP71} and Bunch and
Kaufman~\cite{BunK77}, while Duff et al.~\cite{DufGRST91} and Fourer
and Mehrotra~\cite{FouM92} have described sparse variants. These
algorithms differ in their selection criteria for the $1 \times 1$ and
$2 \times 2$ pivot blocks. In our case, the presence of the diagonal
elements of size $\Theta(\mu^{-1})$ (from the submatrix $D_{\cN} =
\LN^{-1} \SN$) and their place in these pivot blocks are crucial to the
result.

We start by stating a general result of Higham~\cite{Hig97} concerning
backward stability that applies to all diagonal pivoting schemes.  We
then examine the Bunch-Kaufman scheme, showing that the large
diagonals appear only as $1 \times 1$ pivots and that this algorithm
satisfies Condition~\ref{cond.1}. (In \cite[Theorem~4.2]{Hig97},
Higham actually proves that the Bunch-Kaufman scheme is backward
stable in the normwise sense, but this result is not applicable to our
context, for the reasons mentioned above.)

Next, we briefly examine the Bunch-Parlett method, showing that it
starts out by selecting all the large diagonal elements in turn as $1
\times 1$ pivots, before going on to factor the remaining matrix,
whose elements are all $O(1)$ in size. This method also satisfies
Condition~\ref{cond.1}. We then examine the sparse diagonal pivoting
approaches of Duff et al.~\cite{DufGRST91} and Fourer and
Mehrotra~\cite{FouM92}, which may not satisfy Condition~\ref{cond.1},
because of the possible presence of $2 \times 2$ pivots in which one
of the diagonals has size $\Theta(\mu^{-1})$. These algorithms can be
modified in simple ways to overcome this difficulty, possibly at the
expense of higher density in the $L$ factor. We then mention Gaussian
elimination with pivoting and refer to previous results in the
literature to show that this approach satisfies
Condition~\ref{cond.1}.
Finally, we state a result like Theorem~\ref{th:cond2} about
convergence of a finite-precision implementation of an algorithm based
on the augmented system form.

Higham~\cite[Theorem~4.1]{Hig97} proves the following result.
\begin{theorem} \label{th:higham}
Let $T$ be an $\bar{n} \times \bar{n}$ symmetric matrix, and let
$\hat{x}$ be the computed solution to the linear system $Tx=b$
produced by a method that yields a factorization of the form
\eqnok{aug:dpiv}, with any diagonal pivoting strategy. Assume that,
during recovery of the solution, the subsystems that involve the $2
\times 2$ diagonal blocks are solved via Gaussian elimination with
partial pivoting. Then we have that
\beq \label{aug:nh}
(T + \Delta T) \hat{x} = b, \sgap
| \Delta T | \le \dbu ( | T | + 
P^T |\hat{L}| | \hat{Y} | | \hat{L}^T | P) + \dbu^2,
\eeq
where $\hat{L}$ and $\hat{Y}$ are the computed factors, and $|A |$
denotes the matrix formed from $A$ by replacing all its elements by
their absolute values.
\end{theorem}

In Higham's result, the coefficient of $\bu$ in the $\dbu$ term is
actually a linear polynomial in the dimension of the system. The
partial pivoting strategy for the $2 \times 2$ systems can actually be
replaced by any method for which the computed solution of $Ry=d$
satisfies $(R + \Delta R) \hat{y} = d$, where $R$ is the $2 \times 2$
matrix in question and $| \Delta R | \le \dbu |R|$. This property was
also key in an earlier paper of S.~Wright~\cite{Wri94a_rev}, who
derived a result similar to Theorem~\ref{th:higham} in the context of
the augmented systems that arise from interior-point methods for
linear programming.

All the procedures below have the property that the growth in the
maximum element size in the remaining submatrix is bounded by a modest
quantity at each individual step of the factorization. (In the case of
Bunch-Kaufman and Bunch-Parlett, this bound averages about 2.6 per
elimination step; see Golub and Van
Loan~\cite[Section~4.4.4]{GolV96}.) As with Gaussian elimination with
partial pivoting, exponential element growth is possible, so that $L$
and $Y$ in \eqnok{aug:dpiv} contain much larger elements than the
original matrix $T$. Such behavior is, however, quite rare and is
confined to pathological cases and certain special problem classes.
In our analysis below, we make the safe assumption that
catastrophic growth of this kind does not occur.

\subsection{The Bunch-Kaufman Procedure} \label{subs:BK}

At each iteration, the Bunch-\\
Kaufman procedure chooses either a $1
\times 1$ or $2 \times 2$ pivot by examining at most two columns of
the remaining matrix, that is, the part of the matrix that remains to
be factored at this stage of the process. It makes use of quantities
$\chi_i$ defined by
\[
\chi_i = \max_{j \, | \, j \neq i} \, | T_{ij} |,
\]
where in this case $T$ denotes the remaining matrix.  We define the
pivot selection strategy for the first step of the factorization
process. The entire algorithm is obtained by applying this procedure
recursively to the remaining submatrix.

\medskip

\btab
\> set $\nu = (1+\sqrt{17})/8$; \\
\> calculate $\chi_1$, and store  
     the index $r$ for which $\chi_1 = |T_{r1}|$; \\
\> {\bf if} $|T_{11}| \ge \nu \chi_1$ \\
\>\> choose $T_{11}$ as a $1 \times 1$ pivot; \\
\> {\bf else} \\
\>\> calculate  $\chi_r$; \\
\>\> {\bf if} $\chi_r | T_{11} | \ge \nu \chi_1^2$ \\
\>\>\> choose $T_{11}$ as a $1 \times 1$ pivot; \\
\>\> {\bf else} {\bf if}  $|T_{rr} | \ge \nu \chi_r$ \\
\>\>\> choose $T_{rr}$ as a $1 \times 1$ pivot; \\
\>\> {\bf else} \\
\>\>\> choose a $2 \times 2$ pivot with diagonals $T_{11}$ and $T_{rr}$; \\
\>\> {\bf end if}\\
\> {\bf end if.}
\etab

\medskip

For each choice of pivot, the permutation matrix $P_1$ is chosen so
that the desired $1 \times 1$ or $2 \times 2$ pivot is in the upper
left of the matrix $P_1 T P_1^T$. If one writes
\[
P_1 T P_1^T = \bmat{cc} R & C^T \\ C & \hat{T} \emat,
\]
where $R$ is the chosen pivot, the first step of the factorization
yields
\beq \label{aug.upd}
P_1 T P_1^T = \bmat{cc} I & \\ CR^{-1} & I \emat 
\bmat{cc} R & \\ & \bar{T} \emat 
\bmat{cc} I & R^{-1} C^T \\ & I \emat,
\eeq
where $\bar{T} = \hat{T} - CR^{-1} C^T$ is the matrix remaining after this
stage of the factorization. 

At the first step of the factorization, the quantities $\chi_1$ and
$\chi_r$ (if calculated) both have size $O(1)$, since the large
elements of this matrix occur only on the diagonal. Since a $2 \times
2$ pivot is chosen only if
\[
|T_{11}| < \nu \chi_1 \sgap \makebox{\rm and} \sgap  |T_{rr} | < \nu \chi_r,
\]
it follows immediately that both diagonals in a $2 \times 2$ pivot
must be $O(1)$. Hence, the pivot chosen by this procedure is one of
three types:
\begin{subequations} \label{aug.pivots}
\beqa
\label{aug.pivot.1}
&& \makebox{\rm  $1 \times 1$ pivot of size $O(1)$;} \\
\label{aug.pivot.2}
&& \makebox{\rm  $2 \times 2$ pivot in which both diagonals
have size $O(1)$;} \\
\label{aug.pivot.3}
&& \makebox{\rm $1 \times 1$ pivot of size $\Theta(\mu^{-1})$.}
\eeqa
\end{subequations}

In fact, the pivots are one of the types \eqnok{aug.pivots} at {\em
all} stages of the factorization, not just the first stage.  The
reason is that the updated matrix $\bar{T}$ in \eqnok{aug.upd} has the
same essential form as the original matrix $T$---its elements are all
of size $O(1)$ except for some large diagonal elements of size
$\Theta(\mu^{-1})$. We demonstrate this claim by showing that the
update $CR^{-1} C^T$ that is applied to the remaining matrix in
\eqnok{aug.upd} is a matrix whose elements are of size at most $O(1)$,
regardless of the type of pivot, so that it does not disturb the
essential structure of the remaining matrix.  When the pivots are of
type \eqnok{aug.pivot.1} and \eqnok{aug.pivot.2}, the standard
argument of Bunch and Kaufman~\cite{BunK77} can be applied to show
that the norm of $ C R^{-1} C^T$ is at most a modest multiple of
$\|C\|$. We know that $\| C \| = O(1)$, since $C$ consists only of
off-diagonal elements, so we conclude that $\| C R^{-1} C^T \| = O(1)$
in this case as claimed.  For the other pivot type
\eqnok{aug.pivot.3}, we have $R = \Theta(\mu^{-1})$ and $C = O(1)$, so
the elements of $CR^{-1} C^T$ have size $O(\mu)$, and the claim holds
in this case too.

In the rest of this subsection, we show by using
Theorem~\ref{th:higham} that Condition~\ref{cond.1} holds for the
Bunch-Kaufman algorithm. In fact, we prove a stronger result: When $T$
in Theorem~\ref{th:higham} is the matrix \eqnok{gen:aug}, the
perturbation matrix $\Delta T$ contains elements of size $\dbu$,
except in those diagonal locations corresponding to the elements of
$D_{\cN}$, where they may be as large as $\dbu/\mu$. Given the bound
on $|\Delta T|$ in \eqnok{aug:nh}, we need only to show that $P^T
|\hat{L}| |\hat{Y}| |\hat{L}|^T P$ has the desired structure. In fact,
it suffices to show that the exact factor product $P^T |{L}| |{Y}|
|{L}|^T P$ has the structure in question, since the difference between
these two products is just $\dbu$ in size.

We demonstrate this claim inductively, using a refined version of the
arguments from Higham~\cite[Section~4.3]{Hig97} for some key points,
and omitting some details. For simplicity, and without loss of
generality, we assume that $P=I$.

When $\bar{n}=1$ (that is, $T$ is $1 \times 1$), we have that ${L}=1$
and $Y=T$, and the result holds trivially. When $\bar{n}=2$, there are
three cases to consider. If the matrix contains no elements of size
$\Theta(\mu^{-1})$, then the analysis for general matrices can be used
to show that $|{L}| |{Y}| |{L}|^T = O(1)$, as required. If either or
both diagonal elements have size $\Theta(\mu^{-1})$, then both pivots
are $1 \times 1$, and the factors have the form
\beq \label{aug:ldform}
L = \bmat{cc} 1 & 0 \\ T_{21}/T_{11} & 1 \emat, \gap
Y = \bmat{cc}  T_{11} & 0 \\ 0 & T_{22} - T_{21}^2/T_{11} \emat.
\eeq
Two cases arise.
\begin{itemize}
\item[(i)] A diagonal of size $O(1)$ is accepted as the first pivot
and moved (if necessary) to the $(1,1)$ position. We then have
\[
|T_{11}| \ge \nu \chi_1 = \nu \chi_2 = \nu | T_{21}|,
\]
and therefore $|T_{21}/T_{11}| \le 1/\nu$ and hence $|T_{21}^2/T_{11}|
\le |T_{21}|/\nu = O(1)$. If the $(2,2)$ diagonal is also $O(1)$, we
have that $L=O(1)$ and $Y=O(1)$, and we are done. Otherwise, $T_{22} =
\Theta(\mu^{-1})$, and so the $(2,2)$ element of $Y$ satisfies this same
estimate. We conclude from \eqnok{aug:ldform} that $|L||Y||L|^T$
also has an $\Theta(\mu^{-1})$ element in the $(2,2)$ position and
$O(1)$ elements elsewhere.

\item[(ii)] A diagonal of size $\Theta(\mu^{-1})$ is accepted as the
first pivot and moved (if necessary) to the $(1,1)$ position. We then
have
\[
|T_{21}/T_{11}| = O(\mu), \gap 
|T_{21}^2/T_{11}| = O(\mu).
\]
It follows from \eqnok{aug:ldform} that
\[
|L||Y||L|^T = 
\bmat{cc} |T_{11}| & |T_{21}| \\ 
|T_{21}| & |T_{22}| + O(\mu) \emat,
\]
which obviously has the desired structure.
\end{itemize}

We now assume that our claim holds for some dimension $\bar{n} \ge 2$,
and we prove that it continues to hold for dimension $\bar{n}+1$. Using
the notation of \eqnok{aug.upd} (assuming that $P_1=I$), and denoting
the factorization of the Schur complement $\bar{T}$ in \eqnok{aug.upd}
by $\bar{T} = \bar{L} \bar{Y} \bar{L}^T$, we have that
\beq \label{aug.upd2}
T = L Y L^T = \bmat{cc} I & \\ CR^{-1} & \bar{L} \emat 
\bmat{cc} R & \\ & \bar{Y} \emat  
\bmat{cc} I & R^{-1} C^T \\ & \bar{L}^T \emat.
\eeq
It follows that
\beq \label{aug.upd3}
|L||Y||L|^T = \bmat{cc} |R| & |R| |R^{-1} C^T | \\
| CR^{-1}| |R| & 
| CR^{-1}| |R|  |R^{-1} C^T |  + |\bar{L}| |\bar{Y}| |\bar{L}|^T \emat.
\eeq
Since, as we mentioned above, the norm of $CR^{-1}C^T$ is at most
$O(1)$, the Schur complement $\bar{T} = \hat{T} - C R^{-1}C^T$ has
size $O(1)$ except for large $\Theta(\mu^{-1})$ elements in the same
locations as in the original matrix. Hence, by our inductive
hypothesis, $|\bar{L}| |\bar{Y}| |\bar{L}|^T$ has a similar structure,
and we need to show only that the effects of the first step of the
factorization \eqnok{aug.upd} do not disturb the desired structure.

For the case in which $R$ is a pivot of type either
\eqnok{aug.pivot.1} and \eqnok{aug.pivot.2},
Higham~\cite[Section~4.3]{Hig97} shows all elements of both $|
CR^{-1}| |R|$ and $| CR^{-1}| |R| |R^{-1} C^T |$ are bounded by a
modest multiple of either $\chi_1$ (if $T_{11}$ was selected as the
pivot because it passed the test $|T_{11}| \ge \nu \chi_1$) or
$(\chi_1 + \chi_r)$, where $r$ is the ``other'' column considered
during the selection process.  In our case, this observation 
implies that both $| CR^{-1}| |R|$ and $|
CR^{-1}| |R| |R^{-1} C^T |$ have size $O(1)$.  By combining these
observations with those of the preceding paragraph, we conclude that
for pivots of types \eqnok{aug.pivot.1} and \eqnok{aug.pivot.2},
``large'' elements of the matrix in \eqnok{aug.upd3} occur only in the
diagonal locations originally occupied by $D_{\cN}$.

For the remaining case---pivots of type \eqnok{aug.pivot.3}---we have
that $C$ has size $O(1)$ while $R^{-1}$ has size $O(\mu)$.  Therefore,
$| CR^{-1}| |R|$ has size $O(1)$ and
$| CR^{-1}| |R| |R^{-1} C^T |$ has size
$O(\mu)$, while $|R|$, which occupies the $(1,1)$ position in the
matrix \eqnok{aug.upd3}, just as it did in the original matrix $T$,
has size $\Theta(\mu^{-1})$. We conclude that the desired structure
holds in this case as well.

We conclude from this discussion that Condition~\ref{cond.1} holds for
the Bunch-Kaufman procedure. We show later that the perturbations
arising from other sources, namely, roundoff and cancellation in the
evaluation of the matrix and right-hand side, also satisfy the
conditions of Corollary~\ref{co:errp}, so this result can be used to
bound the error in the computed steps.

Finally, we note that it is quite possible for pivots of types
\eqnok{aug.pivot.1} and \eqnok{aug.pivot.2} to be chosen while
diagonal elements of size $\Theta(\mu^{-1})$ still remain in the
submatrix. Therefore, a key assumption of the analysis of Forsgren,
Gill, and Shinnerl~\cite[Theorem~4.4]{ForGS94}---namely, that all the
diagonals of size $\Theta(\mu^{-1})$ are chosen as $1 \times 1$ pivots
before any of the other diagonals are chosen---may not be satisfied by
the Bunch-Kaufman procedure.

\subsection{The Bunch-Parlett Procedure} \label{subs:BP}

The Bunch-Parlett procedure is conceptually simpler but more expensive
to implement than Bunch-Kaufman, since it requires $O(n^2)$ (rather
than $O(n)$) comparisons at each step of the factorization.  The pivot
selection strategy is as follows.

\medskip

\btab
\> set $\nu = (1+\sqrt{17})/8$; \\
\> calculate $\chi_{\rm off} = | T_{rs} | = \max_{i \neq j} | T_{ij}|$,
$\chi_{\rm diag} = |T_{pp}| = \max_i |T_{ii}|$; \\
\> {\bf if} $\chi_{\rm diag} \ge \nu \chi_{\rm off}$ \\
\>\> choose $T_{pp}$ as the $1 \times 1$ pivot; \\
\> {\bf else} \\
\>\> choose the $2 \times 2$ pivot whose off-diagonal element is $T_{rs}$; \\
\> {\bf end if.} 
\etab

\medskip

The elimination procedure then follows as in \eqnok{aug.upd}.

It is easy to show that the Bunch-Parlett procedure starts by
selecting all the diagonals of size $\Theta(\mu^{-1})$ in turn as $1
\times 1$ pivots. (Because of this property, it satisfies the key
assumption of \cite{ForGS94} mentioned at the end of the preceding
section.) The update $CR^{-1} C^T$ generated by each of these pivot
steps has size only $O(\mu)$, so the matrix that remains after this
phase of the factorization contains only $O(1)$ elements. The
remaining pivots are then a combination of types \eqnok{aug.pivot.1}
and \eqnok{aug.pivot.2}.

By using the arguments of the preceding subsection in a slightly
simplified form, we can show that Condition~\ref{cond.1} holds for
this procedure as well.


\subsection{Sparse Diagonal Pivoting} \label{subs:Duff}

For large instances of \eqnok{nlp}, the Bunch-\\
Kaufman and
Bunch-Parlett procedures are usually inefficient because they do not
try to maintain sparsity in the lower triangular factor $L$. Sparse
variants of these algorithms, such as those of Duff et
al.~\cite{DufGRST91} and Fourer and Mehrotra~\cite{FouM92}, use pivot
selection strategies that combine stability considerations with
Markowitz-like estimates of the amount of fill-in that a candidate
pivot will cause in the remaining matrix.

At each stage of the factorization, both algorithms examine a roster
of possible $1 \times 1$ and $2 \times 2$ pivots, starting with those
that would create the least fill-in. As soon as a pivot is found that
meets the stability criteria described below, it is accepted. Both
algorithms prefer to use $1 \times 1$ pivots where possible.

For candidate $1 \times 1$ pivots, Duff et
al.~\cite[p.~190]{DufGRST91} use the following stability criterion:
\beq \label{1x1}
|R^{-1}| \| C \|_{\infty} \le \rho,
\eeq
where the notation $R$ and $C$ is from \eqnok{aug.upd} and $\rho \in
[2,\infty)$ is some user-selected parameter that represents the
tolerable growth factor at each stage of the factorization.  For a $2
\times 2$ pivot, the criterion is
\beq \label{2x2}
|R^{-1}| \bmat{c} \| C_{\cdot,1} \|_{\infty} \\  
\| C_{\cdot,2} \|_{\infty} \emat \le \bmat{c} \rho \\ \rho \emat,
\eeq
where $C_{\cdot,1}$ and $C_{\cdot,2}$ are the two columns of $C$. The
stability criteria of Fourer and Mehrotra~\cite{FouM92} are similar.

As they stand, the stability tests \eqnok{1x1} and \eqnok{2x2} do not
necessarily
restrict the choice of pivots to the three types
\eqnok{aug.pivots}. If a $1 \times 1$ pivot of size $\Theta(\mu^{-1})$
is ever considered for structural reasons, it will pass the test
\eqnok{1x1} (the left-hand side of this expression will have size
$O(\mu)$) and therefore will be accepted as a pivot. However, it is
possible that $2 \times 2$  pivots in which one or both diagonals have size
$\Theta(\mu^{-1})$ may pass the test \eqnok{2x2} and may therefore be
accepted. Although the test \eqnok{2x2} ensures that the size of the
update $CR^{-1} C^T$ is modest (so that the update $\bar{T} = \hat{T}
- CR^{-1} C^T$ does not disturb the large-diagonal structure of
$\hat{T}$), there is no obvious assurance that the matrix
$|L||Y||L|^T$ in \eqnok{aug.upd3} mirrors the structure of $|T|$, in
terms of having the large diagonal elements in the same locations. The
terms $|CR^{-1}||R|$ and $|CR^{-1}||R||R^{-1} C^T|$ in
\eqnok{aug.upd3} may not have size $O(1)$, as they do for pivots of
the three types \eqnok{aug.pivots} arising from the Bunch-Kaufman and
Bunch-Parlett selection procedures.

The Fourer-Mehrotra algorithm does, however, rule out the possibility
of a $2 \times 2$ pivot in which {\em both} diagonals are of size
$\Theta(\mu^{-1})$. It considers a $2 \times 2$ candidate only if one
of its diagonal elements has previously been considered as a $1 \times
1$ pivot but failed the stability test. However, if either of the
diagonals had been subjected to the test \eqnok{1x1}, they would have
been accepted, as noted in the preceding paragraph, so this situation
cannot occur.

If the sparse algorithms are modified to ensure that all pivots have
one of the three types \eqnok{aug.pivots}, and all continue to satisfy
the stability tests \eqnok{1x1} or \eqnok{2x2}, then simple arguments
(simpler than those of Section~\ref{subs:BK}!) can be applied to show
that Condition~\ref{cond.1} is satisfied. One possible modification
that achieves the desired efffect is to require that a $2 \times 2$
pivot be allowed only if {\em both} its diagonals have previously been
considered as $1 \times 1$ pivots but failed the stability test
\eqnok{1x1}.

\subsection{Gaussian Elimination} \label{subs:GE}

Another possibility for solving the system \eqnok{gen:aug} is to
ignore its symmetry and apply a Gaussian elimination algorithm, with
row and/or column pivoting to preserve sparsity and prevent excessive
element growth. Such a strategy satisfies Condition~\ref{cond.1}. In
\cite{Wri93c}, the author uses a result of Higham~\cite{Hig96} to show
that the effects of the large diagonal elements are essentially
confined to the columns in which they appear. Assuming that the pivot
sequence is chosen to prevent excessive element growth in the
remaining matrix, and using the notation of \eqnok{gen.30} and
\eqnok{esizes}, we can account for the effects of roundoff error in
Gaussian elimination with perturbations in the coefficient matrix that
satisfy the following estimates:
\[
E_{11}, E_{21}, E_{31}, E_{12}, E_{22}, E_{32} = \dbu, \sgap
E_{13}, E_{23}, E_{33} = \dbu/\mu.
\]
These certainly satisfy the conditions \eqnok{esizes}, so
Condition~\ref{cond.1} holds.

\subsection{Local Convergence with the Computed Steps} \label{sec:aug.conv}

We can now state a formal result to show that when the evaluation
errors are taken into account as well as the roundoff errors from the
factorization/solve procedure discussed above, the accuracies of the
computed steps obtained from the procedure {\bf augmented},
implemented in finite precision, satisfy the same estimates as for the
corresponding steps obtained from the procedure {\bf condensed}. The
result is analogous to Theorem~\ref{th:cond1}.
\begin{theorem} \label{th:cond1.aug}
Suppose that Assumption~\ref{ass:base} holds.  Then when the step \\
$(\wDz, \wDl, \wDs)$ is calculated in a finite-precision environment
by using the procedure {\bf augmented}, where the algorithm used to
solve \eqnok{gen:aug} satisfies Condition~\ref{cond.1}, we have
\begin{subequations} \label{aug.p16}
\beqa
\label{aug.p16a}
(\Dz - \wDz, U^T (\Dl_{\cB} - \wDl_{\cB}), \Ds - \wDs) &=&
\dbu, \\
\label{aug.p16b}
V^T (\Dl_{\cB} - \wDl_{\cB}) &=& \dbu / \mu, \\
\label{aug.p16c}
\Dl_{\cN} - \wDl_{\cN} &=& \dbu.
\eeqa
\end{subequations}
\end{theorem}
\begin{proof}
The proof follows from Corollary~\ref{co:errp} when we show that the
perturbations to \eqnok{gen:aug} from all sources---evaluation of the
matrix and right-hand side as well as the factorization/solution
procedure---satisfy the bounds required by this earlier result.

Because of Condition~\ref{cond.1}, perturbations arising from the
factorization/solution procedure satisfy the bounds \eqnok{esizes}.
The expressions \eqnok{cancel} show that the errors arising from
evaluation of $\cL_{zz} (z,\lambda)$, $\cL_z(z,\lambda)$, $\nablag(z)$,
and $g(z)$ are all of size $\dbu$, and hence they too satisfy the
required bounds. Similarly to \eqnok{cancelD}, evaluation of $D_{\cB}$
and $D_{\cN}$ yields errors of relative size $\dbu$, that is,
\begin{subequations} \label{cancelD.1}
\beqa \label{cancelD.B.1}
\makebox{\rm computed } D_{\cB}  & \leftarrow & D_{\cB} + G_{\cB},
 \sgap G_{\cB} = \mu \dbu, \\
\label{cancelD.N.1}
\makebox{\rm computed } D_{\cN} & \leftarrow & D_{\cN} + G_{\cN}, 
\sgap G_{\cN} = \dbu/\mu,
\eeqa
\end{subequations}
where $G_{\cB}$ and $G_{\cN}$ are diagonal matrices. 

We now obtain all the estimates in \eqnok{aug.p16} by a direct
application of Corollary~\ref{co:errp}, with the exception of the
estimate for $(\Ds - \wDs)$. Since the expressions for recovering
$\Ds$ are identical in procedures {\bf condensed} and {\bf augmented},
we can apply expression \eqnok{cond.p14} from
Section~\ref{sec:cond:quant} to deduce that the desired estimate holds
for this component as well.
\end{proof}

The only difference between the error estimates of
Theorem~\ref{th:cond1} for the condensed system and those obtained
above for the augmented system is that the $\wDl_{\cN}$ components are
slightly less accurate in the augmented case. If we work through the
analysis of Section~\ref{sec:cond:local} with the estimate
\eqnok{aug.p16c} replacing \eqnok{cond.p16c}, we find that the main
results are unaffected. Therefore, we conclude this section by stating
without proof a result similar to Theorem~\ref{th:cond2}.
\begin{theorem} \label{th:aug2}
Suppose that all the assumptions of Theorem~\ref{th:cond2} hold,
except that the step $(\wDz, \wDl, \wDs)$ is calculated by using the
procedure {\bf augmented} with a factorization/solution algorithm that
satisfies Condition~\ref{cond.1}. Then the conclusions of
Theorem~\ref{th:cond2} hold.
\end{theorem}

\section{Numerical Illustration}
\label{sec:numerics}

We illustrate the results of Sections~\ref{sec:pdip.cond} and
\ref{sec:pdip.aug} using the two-variable example
\eqnok{numex}. Consider a simple algorithm that takes steps satisfying
\eqnok{pdip.orig} with $t$ set rather arbitrarily to $t=\mu^2 e$. (The
search directions thus used are like those generated in the later
stages of a practical primal-dual algorithm such as Mehrotra's
algorithm~\cite{Meh92a}.) We start this algorithm from the point
\[
z_0=(1/30,1/9)^T, \sgap \lambda_0 = (1,1/5)^T, \sgap s_0=(1/10, 1/2).
\]
(It is easy to check that the conditions \eqnok{pdip.central} are
satisfied at this point for a modest value of $C$.) At each step we
calculated $\walfmax$, defined in Section~\ref{sec:cond:local}, and
took an actual step of $.99 \hat{\alpha}_{\rm max}$.

We programmed the method in Matlab, using double-precision arithmetic.
In our first experiment, we solved the formulation \eqnok{gen:aug} of
the linear equations using Matlab's standard Gaussian elimination
solver for general systems of linear equations, which was analyzed in
Section~\ref{subs:GE}. From Theorem~\ref{th:cond1.aug}, the estimates
\eqnok{aug.p16} apply to this case.

\begin{table}
{\centering%
\begin{tabular}{|cr|ccc|rc|} \hline
\mbox{iter} & $\log \mu$ & $\log \| \wDz \|$ & $\log \| U^T \wDl_{\cB} \|$ &
$\log \| V^T \wDl_{\cB} \|$ &  $ \walfmax$ & $\lambda^T$  \\ \hline
0 & -1.0 & -0.9 & -1.9 & -1.9 & .9227 & (1.00,.20) \\
1 & -2.7 & -1.5 & -1.3 & -1.2 & .9193 & (0.99,.19) \\
\vdots &&&&&& \\
5 & -9.4 & -6.7 & -6.3 & -4.6 & 1.0 & (1.04,.23) \\
6 & -11.4 & -8.7 & -8.3 & -5.9 & 1.0 & (1.04, .23) \\
7 & -13.4 & -10.7 & -10.3 & -3.8 & .9999 & (1.04,.23) \\
8 & -15.4 & -12.7 & -12.3 & -1.2 & .9439 & (1.04,.23) \\
9 & -17.1 & -13.9 & -13.4 & -0.6 & .9723 & (1.10,.20) \\ \hline
\end{tabular}
}
\caption{Details of iteration sequence for PDIP applied to (\protect\ref{numex}), with steps computed by solving the augmented system.\label{ta:pdip.aug1}}
\end{table}

Results are tabulated in Table~\ref{ta:pdip.aug1}. Note first the size
of the component $\| V^T \wDl_{\cB} \|$, which grows as $\mu$
decreases below $\bu^{1/2}$, in accordance with 
\eqnok{aug.p16b}. (We cannot tabulate the difference $\| V^T
(\wDl_{\cB} - \Dl_{\cB} \|)$ because of course we do not know the true
step $(\Dz, \Dl, \Ds)$, but since the true step has size $O(\mu)$
(Corollary~\ref{co:exact}), the error is dominated by the term $V^T
\wDl_{\cB}$ in any case.) As predicted by \eqnok{cond.p27}, the
maximum step $\walfmax$ becomes significantly smaller than $1$ as
$\mu$ is decreased below $\bu^{1/2}$. As indicated by
\eqnok{cond.p28}, however, good progress still can be made along this
direction (in the sense of reducing $\mu$ and the norms of the
residuals $r_f$ and $r_g$) almost until $\mu$ reaches the level of
$\bu$. In fact, between iterations 5 and 8 we see the reduction factor
of $100$ that we would expect by moving a distance of $.99$ along a
direction that is close to the pure Newton direction.  The component
with the large error---$ V^T \wDl_{\cB}$---does not interfere
significantly with rapid convergence, but only causes the $\lambda$
iterates to move tangentially to $\cS_{\lambda}$. This effect may be
noted in the final iterate where the value of $\lambda$ changes
significantly. In some cases, however, when the current $\lambda$ is
near the edge of the set $\cS_{\lambda}$, this error may result in a
severe curtailment of the step length.

\begin{table}
{\centering%
\begin{tabular}{|cr|ccc|rc|} \hline
\mbox{iter} & $\log \mu$ & $\log \| \wDz \|$ & $\log \| U^T \wDl_{\cB} \|$ &
$\log \| V^T \wDl_{\cB} \|$ &  $ \walfmax$ & $\lambda^T$  \\ \hline
0 & -1.0 & -0.9 & -1.9 & -1.9 & .9227 & (1.00,.20) \\
1 & -2.7 & -1.5 & -1.3 & -1.2 & .9193 & (0.99,.19) \\
\vdots &&&&&& \\
5 & -9.4 & -6.7 & -6.3 & -4.6 & 1.0 & (1.04,.23) \\
6 & -11.4 & -8.7 & -8.3 & -5.7 & 1.0 & (1.04, .23) \\
7 & -13.4 & -10.7 & -10.3 & -8.3 & 1.0 & (1.04,.23) \\
8 & -15.4 & -12.7 & -12.4 & -10.3 & 1.0 & (1.04,.23) \\ 
9 & -17.4 & -14.7 & -13.3 & -12.3 & 1.0 & (1.04,.23) \\ 
\hline
\end{tabular}
}
\caption{Details of iteration sequence for PDIP applied to (\protect\ref{numex}), with steps computed by solving the condensed system.\label{ta:pdip.con1}}
\end{table}

Next, we performed the same experiment using the condensed formulation
\eqnok{pdip.cond} of the linear system, as described in
Section~\ref{sec:pdip.cond}. Results are shown in
Table~\ref{ta:pdip.con1}. The main difference with
Table~\ref{ta:pdip.aug1} is that there is no increase in the value $\|
V^T \wDl_{\cB} \|$ as $\mu$ approaches unit roundoff; this component
appears to decrease at the same rate as the other step components.
This observation can be explained by our analysis of the case in which
the cancellation error term $f_{\cB}$ incurred in the evaluation of
$\gB(z)$ satisfies \eqnok{improved.vtf}.  We calculated the product
$V^T (\gB(z) + f_{\cB})$ (the product of $V$ with our computed version
of $\gB(z)$) and found it to be exactly zero on iterations 7, 8, and
9. Therefore, using Taylor's theorem, \eqnok{gen.3}, and
Theorem~\ref{th:est3}, we have
\[
V^T f_{\cB} = - V^T \gB(z) = -V^T \nablagB(z^*) (z-z^*)  + O(\| z-z^*\|^2) =
O(\mu^2).
\]
Hence, \eqnok{cond.p16b.improved} together with
Corollary~\ref{co:exact} shows that $V^T \wDl_{\cB} = O(\mu)$, which
is consistent with the results in Table~\ref{ta:pdip.con1}.  Note too
that because of the higher accuracy in the $V^T \wDl_{\cB}$ component,
the maximum step length stays very close to $1$ during the last few
iterations.  By comparing Tables~\ref{ta:pdip.aug1} and
\ref{ta:pdip.con1}, however, we can verify that the convergence of
$\mu$ to zero, and of the iterates to the solution set, is not
materially affected by the presence or absence of the large error in
$V^T \wDl_{\cB}$.

To show that the lack of cancellation effects in
Table~\ref{ta:pdip.con1} cannot be assumed in general, we modified
problem \eqnok{numex} slightly, changing the second constraint to
\beq \label{numex.modified}
g_2(z) \defeq 
\frac{2}{3 \sqrt{5}} (z_1- \sqrt{5})^2 + z_2^2 - \frac{2 \sqrt{5}}{3} \le 0.
\eeq
The primal and dual solutions remain unchanged, and we ran the
condensed-equations version of the algorithm from the same
starting point as above. Results are shown in Table~\ref{ta:pdip.con2}. We
observed that $\gB(z)$ did not escape cancellation errors in this
instance and, as in Table~\ref{ta:pdip.aug1}, we observe significant
errors in $V^T \wDl_{\cB}$ that do not materially affect the
convergence of the algorithm to the solution set.

\begin{table}
{\centering%
\begin{tabular}{|cr|ccc|rc|} \hline
\mbox{iter} & $\log \mu$ & $\log \| \wDz \|$ & $\log \| U^T \wDl_{\cB} \|$ &
$\log \| V^T \wDl_{\cB} \|$ &  $ \walfmax$ & $\lambda^T$  \\ \hline
0 & -1.0 & -0.9 & -2.1 & -2.3 & .9161 & (1.00,.20) \\
1 & -2.7 & -1.5 & -1.3 & -1.4 & .8872 & (0.99,.20) \\
\vdots &&&&&& \\
5 & -7.6 & -5.7 & -5.7 & -4.2 & .9999 & (.93,.29) \\
6 & -9.5 & -7.7 & -7.7 & -6.3 & 1.0 & (.93, .29) \\
7 & -11.5 & -9.7 & -9.7 & -4.3 & .9999 & (.93, .29) \\
8 & -13.5 & -11.7 & -11.5 & -2.6 & .9960 & (.93,.29) \\ 
9 & -15.3 & -13.5 & -11.7 & -0.6 & .7386 & (.93,.29) \\ 
\hline
\end{tabular}
}
\caption{Details of iteration sequence for PDIP applied to (\protect\ref{numex}), (\protect\ref{numex.modified}), with steps computed from the condensed
system.\label{ta:pdip.con2}}
\end{table}

\section{Summary and Conclusions} \label{sec:conclusions}

In this paper, we have analyzed the finite-precision implementation of
a primal-dual interior point method whose convergence rate is
theoretically superlinear.  We have made the standard assumptions that
appear in most analyses of local convergence of nonlinear programming
algorithms and path-following algorithms, with one significant
exception: The assumption of linearly independent active constraint
gradients is replaced by the weaker Mangasarian-Fromovitz constraint
qualification, which is equivalent to boundedness of the set of
optimal Lagrange multipliers. Because of this assumption, it is
possible that all reasonable formulations of the step equations---the
linear system that needs to be solved to obtain the search
direction---are ill conditioned, so it is not obvious that the
numerical errors that occur when this system is solved in finite
precision do not eventually render the computed search direction
useless. We show that although the error in the computed step may
indeed become large as $\mu$ decreases, most of the error is
restricted to a subspace that does not matter, namely, the null space
of the matrix $\nablagB(z^*)$ of first derivatives of the active
constraints. Although this error causes the computed iterates to
``slip'' in a tangential direction to the optimal Lagrange multiplier
set, it does not interfere with rapid convergence of the iterates to
the primal-dual solution set.

We found that the centrality conditions \eqnok{pdip.central}, which
are usually applied in path-following methods, played a crucial role
in the analysis, since they enabled us to establish the estimates
\eqnok{eq:14} in Lemma~\ref{lem:est2} concerning the sizes of the
basic and nonbasic components of $s$ and $\lambda$ near the solution
set. The analysis of Section~\ref{accuracy}, culminating in
Corollary~\ref{co:errp}, finds bounds on the errors induced in step
components by certain structured perturbations of the step equations.
We show in the same section that the exact step is $O(\mu)$, allowing
the local convergence analysis of Ralph and Wright~\cite{RalW96b} to
be extended from convex programs to nonlinear programs.

In Sections~\ref{sec:pdip.cond} and \ref{sec:pdip.aug} we apply the
general results of Section~\ref{accuracy} to the two most obvious ways
of formulating and solving the step equations; namely, as a
``condensed'' system involving just the primal variables $z$, or as an
``augmented'' system involving both $z$ and the Lagrange multipliers
$\lambda$. In each case, the errors introduced in finite-precision
implementation have the structure of the perturbations analyzed in
Section~\ref{accuracy}, so the error bounds obtained in
Corollary~\ref{co:errp} apply. In Section~\ref{sec:cond:local} (whose
analysis also applies to the computed solutions analyzed in
Section~\ref{sec:pdip.aug}), we show that the potentially large error
component discussed above does not interfere appreciably with the
near-linear decrease of the quantities $\mu$, $r_f$, and $r_g$ to zero
along the computed steps, indicating that until $\mu$ becomes quite
close to $\bu$, the convergence behavior predicted by the analysis of
the ``exact'' algorithm will be observed in the finite-precision
implementation. We conclude in Section~\ref{sec:numerics} with a
numerical illustration of our major observations on a simple problem
with two variables and two constraints, first introduced in
Section~\ref{notation}.

\section*{Acknowledgments}

Many thanks are due to an anonymous referee for close and careful
readings of various versions of the paper and for many helpful
suggestions.

%
%

\appendix

\section*{}

{\bf Justification of the Estimates (\protect\ref{cond.p27}) and (\protect\ref{cond.p28}).}

To prove \eqnok{cond.p27}, we use analysis similar to that of
S.~Wright~\cite{Wri96b}. From the definition \eqnok{def.mu} of
$\mu$, and the centrality condition \eqnok{pdip.central.3}, we have
that
\[
\lambda_i s_i = \Theta(\mu), \sgap \makebox{\rm for all $i=1,2,\dots,m$.}
\]
Hence, from the third block row of \eqnok{pdip.orig} and the
assumption \eqnok{tpert.est} on the size of $\tpert$,
we have that 
\beq \label{cond.p30}
\frac{\Dl_i}{\lambda_i} + \frac{\Ds_i}{s_i} = -1 - \frac{t_i}{s_i \lambda_i}
= -1 + O(\mu), \gap \makebox{\rm for all $i=1,2,\dots,m$.}
\eeq
We have from Lemma~\ref{lem:est2} and \eqnok{est:exact} that $\Dl_i /
\lambda_i = O(\mu)$ for all $i \in \cB$. Hence, by using \eqnok{eq:14a} from
\eqnok{lem:est2} together with \eqnok{cond.p30}, we obtain
\beq \label{cond.p31}
\Ds_i = -s_i + O(\mu^2), \gap \makebox{\rm for all $i \in \cB$.}
\eeq
For the computed step components $\wDs_{\cB}$, we have by combining
\eqnok{cond.p16a} with \eqnok{cond.p31} that
\beq \label{cond.p32}
\wDs_i = -s_i + \dbu + O(\mu^2), \gap \makebox{\rm for all $i \in \cB$.}
\eeq
Therefore, if $s_i + \alpha \wDs_i=0$ for some $i \in \cB$ and some
$\alpha \in [0,1]$, we have by using \eqnok{eq:14a} again that
\beqa
\nonumber
&& s_i + \alpha(-s_i + \dbu + O(\mu^2)) = 0 \\
\nonumber
& \Rightarrow & (1-\alpha) s_i =\dbu + O(\mu^2) \\
\label{cond.p33}
& \Rightarrow & (1-\alpha) = \dbu/\mu + O(\mu), 
\gap \makebox{\rm for any $i \in \cB$.}
\eeqa
Meanwhile, for $i \in \cN$, we have from Lemma~\ref{lem:est2},
\eqnok{est:exact}, and \eqnok{cond.p16a} that
\beq \label{cond.p34}
s_i + \alpha \wDs_i >0, \gap 
\makebox{\rm for all $i \in \cN$ and all $\alpha \in [0,1]$,}
\eeq
so the components $\wDs_{\cN}$ do not place a limit on the step length
bound $\walfmax$. For the components $\wDl_{\cN}$, we have by using
Lemma~\ref{lem:est2}, \eqnok{est:exact}, \eqnok{cond.p16c}, and
\eqnok{cond.p30} that
\[
\wDl_i =  -\lambda_i + \mu \dbu + O(\mu^2), 
\gap \makebox{\rm for all $i \in \cN$.}
\]
Therefore, if $\lambda_i + \alpha \wDl_i=0$ for some $i \in \cN$ and
some $\alpha \in [0,1]$, we have by arguing as in \eqnok{cond.p33}
that 
\beq \label{cond.p35}
1-\alpha = \dbu + O(\mu).
\eeq
Finally, for $i \in \cB$, we have from Lemma~\ref{lem:est2} that $\lambda_i
= \Theta(1)$, while from \eqnok{est:exact}, \eqnok{cond.p16a}, and
\eqnok{cond.p16b}, we have that
\beq \label{cond.p36}
\Dl_i = O(\mu), \gap \wDl_i = O(\mu) + \dbu/\mu, 
\gap \makebox{\rm for all $i \in \cB$.}
\eeq
Therefore, we have for  $\mu \gg \bu$ that 
\beq \label{cond.p37}
\lambda_i + \alpha \wDl_i >0, \gap 
\gap \makebox{\rm for all $i \in \cB$ and all $\alpha \in [0,1]$.}
\eeq
By combining the observations \eqnok{cond.p33}, \eqnok{cond.p34},
\eqnok{cond.p35}, and \eqnok{cond.p37}, we conclude that
there is a value $\walfmax$ satisfying 
\[
\walfmax \in [0,1], \gap 1-\walfmax = \dbu/\mu + O(\mu)
\]
such that
\[
(\lambda,s) + \alpha (\wDl, \wDs) >0, \gap 
\makebox{\rm for all $\alpha \in [0,\walfmax]$,}
\]
proving the claim \eqnok{cond.p27}.  By making various simplifications
to the analysis above, it is easy to show that \eqnok{cond.p25} holds
as well.

We now prove the claims \eqnok{cond.p28} concerning the changes in the
feasibility and duality measures along the computed step. 

From \eqnok{lagrange}, \eqnok{pdip.central.1}, and the first block row
of \eqnok{pdip.orig}, we have
\beqa
\nonumber
\lefteqn{r_f(z+\alpha \wDz, \lambda + \alpha \wDl)} \\
\nonumber
&=& \cL_z (z+\alpha \wDz, \lambda + \alpha \wDl) \\
\nonumber
&=& \cL_z(z,\lambda) + \alpha \cL_{zz} (z,\lambda) \wDz + \alpha
\nablag(z) \wDl + O(\alpha^2 \| \wDz \|^2) \\
\nonumber 
&=& (1-\alpha) \cL_z(z,\lambda) + \alpha \cL_{zz} (z,\lambda) (\wDz - \Dz) 
+ \alpha \nablag_{\cB}(z) (\wDl_{\cB} - \wDl_{\cB}) \\
\label{cond.p40}
&& \gap + \alpha \nablag_{\cN}(z) (\wDl_{\cN} - \Dl_{\cN}) + O( \alpha^2
\| \wDz\|^2).
\eeqa
From \eqnok{est:exact} and \eqnok{cond.p16a}, we have $\wDz = \dbu +
O(\mu)$, so for $\mu \gg \bu$ and $\alpha \in [0,1]$, we have
\beq \label{cond.p41}
\alpha^2 \| \wDz \|^2  = O(\mu^2).
\eeq
From the definition \eqnok{gen.3} of the SVD of $\nablag_{\cB}(z^*)$,
Theorem~\ref{th:est3}, and \eqnok{cond.p16a}, we have that
\beqa
\nonumber
\nablag_{\cB}(z) (\wDl_{\cB} - \Dl_{\cB}) &=& 
\nablag_{\cB}(z^*) (\wDl_{\cB} - \Dl_{\cB}) + 
O(\| z-z^*\| \| \wDl_{\cB} - \Dl_{\cB} \|) \\
\nonumber
&=& \hat{U} \Sigma U^T (\wDl_{\cB} - \Dl_{\cB})  + 
O(\mu) \dbu/\mu  \\ 
\label{mfcq.crap}
&=& \dbu.
\eeqa
Note that the larger error \eqnok{cond.p16b} in the component $V^T
(\wDl_{\cB} - \Dl_{\cB})$, which is present when MFCQ is satisfied but not
when LICQ is satisfied, does not enter into the estimate
\eqnok{mfcq.crap}.
By substituting this estimate into \eqnok{cond.p40} together with
estimates for $\wDz - \Dz$ and $\wDl_{\cN} - \Dl_{\cN}$ from
\eqnok{cond.p16}, we obtain that
\[
r_f(z+\alpha \wDz, \lambda + \alpha \wDl) = (1-\alpha) r_f + \dbu + O(\mu^2),
\]
verifying our claim \eqnok{cond.p28b}. The potentially large error
\eqnok{cond.p16b} does not affect rapid decrease of the $r_f$
component along the computed search direction.

For the second feasibility measure $r_g$, we have from
\eqnok{pdip.central.2}, the second block row of \eqnok{pdip.orig}, and
the estimates \eqnok{cond.p16a} and \eqnok{cond.p41} that
\beqas
\lefteqn{r_g(z+\alpha \wDz, s+\alpha \wDs)} \\
&=& g(z+\alpha \wDz) + s+\alpha \wDs \\
&=& g(z) + \alpha \nablag(z)^T \wDz + s + \alpha \wDs + O(\alpha^2 \| \wDz\|^2) \\
&=& (1-\alpha) (g(z)+s) + \alpha \nablag(z)^T (\wDz - \Dz) +
\alpha (\wDs - \Ds) + O(\mu^2) \\
&=& (1-\alpha) r_g + \dbu + O(\mu^2),
\eeqas
verifying \eqnok{cond.p28c}.

To examine the change in $\mu$, we look at the change in each pairwise
product $\lambda_i s_i$, $i=1,2,\dots,m$. We have
\beqa
\nonumber
\lefteqn{(\lambda_i + \alpha \wDl_i)(s_i + \alpha \wDs_i)} \\
\nonumber 
&=& \lambda_i s_i + \alpha(s_i \wDl_i + \lambda_i \wDs_i) 
+ \alpha^2 \wDs_i \wDl_i \\
\label{cond.p44}
&=& \lambda_i s_i + \alpha(s_i \Dl_i + \lambda_i \Ds_i) + \alpha s_i
(\wDl_i - \Dl_i) + \alpha \lambda_i (\wDs_i - \Ds_i) \\
\nonumber 
&& \gap + \alpha^2 \wDl_i \wDs_i.
\eeqa
From the last block row in \eqnok{pdip.orig}, the estimate
$\tpert=O(\mu^2)$ \eqnok{tpert.est}, and the estimate
\eqnok{est:exact} of the exact step, we have
\beq \label{cond.p45}
 \lambda_i s_i + \alpha(s_i \Dl_i + \lambda_i \Ds_i) =
(1-\alpha)  \lambda_i s_i  + O(\mu^2).
\eeq
From \eqnok{est:exact} and \eqnok{cond.p16}, we have
\beq \label{cond.p46}
\wDl_i \wDs_i = (\dbu / \mu + O(\mu)) (O(\mu) + \dbu)  =
\dbu + O(\mu^2),
\eeq
since $\mu \gg \bu$. For $i \in \cB$, we have from Lemma~\ref{lem:est2},
\eqnok{cond.p16a}, and \eqnok{cond.p16b} that
\beq \label{cond.p47}
s_i(\wDl_i - \Dl_i) = O(\mu) \dbu/\mu  =
\dbu, \sgap \makebox{\rm for all $i \in \cB$}.
\eeq
For $i \in \cN$, we have from Lemma~\ref{lem:est2} and \eqnok{cond.p16c} that
\beq \label{cond.p48}
s_i(\wDl_i - \Dl_i) = \mu \dbu, 
\sgap \makebox{\rm for all $i \in \cN$}.
\eeq
For the remaining term $\lambda_i (\wDs_i - \Ds_i)$, we have from
Lemma~\ref{lem:est2} and \eqnok{cond.p16a} that 
\beq \label{cond.p49}
\lambda_i ( \wDs_i - \Ds_i) = \dbu, \sgap 
\makebox{\rm for all $i=1,2,\dots,m$.}
\eeq
By substituting \eqnok{cond.p45}--\eqnok{cond.p49} into \eqnok{cond.p44}, 
we obtain
\beq \label{cond.p49a}
(\lambda_i + \alpha \wDl_i)(s_i + \alpha \wDs_i) = 
(1-\alpha) \lambda_i s_i  + \dbu + O(\mu^2), \sgap
\makebox{\rm all $i=1,2,\dots,m$.}
\eeq
Therefore, by summing over $i$ and using  \eqnok{def.mu}, we obtain
\eqnok{cond.p28a}.

\bibliography{refs}

\begin{thebibliography}{10}

\bibitem{AndGMX96}
{\sc E.~D. Andersen, J.~Gondzio, C.~{M\'esz\'aros}, and X.~Xu}, {\em
  Implementation of interior-point methods for large scale linear programming},
  in Interior Point Methods in Mathematical Programming, T.~Terlaky, ed.,
  Kluwer Academic Publishers, 1996, ch.~6, pp.~189--252.

\bibitem{BunK77}
{\sc J.~Bunch and L.~Kaufman}, {\em Some stable methods for calculating inertia
  and solving symmetric linear systems}, Mathematics of Computation, 31 (1977),
  pp.~163--179.

\bibitem{BunP71}
{\sc J.~R. Bunch and B.~N. Parlett}, {\em Direct methods for solving symmetric
  indefinite systems of linear equations}, SIAM Journal on Numerical Analysis,
  8 (1971), pp.~639--655.

\bibitem{ByrLN98}
{\sc R.~H. Byrd, G.~Liu, and J.~Nocedal}, {\em On the local behavior of an
  interior-point method for nonlinear programming}, {OTC} Technical Report
  98/02, Optimization Technology Center, January 1998.

\bibitem{PCx99}
{\sc J.~Czyzyk, S.~Mehrotra, M.~Wagner, and S.~J. Wright}, {\em {PCx}: An
  interior-point code for linear programming}, Optimization Methods and
  Software, 11/12 (1999), pp.~397--430.

\bibitem{Deb52}
{\sc G.~Debreu}, {\em Definite and semidefinite quadratic forms}, Econometrica,
  20 (1952), pp.~295--300.

\bibitem{DufGRST91}
{\sc I.~S. Duff, N.~I.~M. Gould, J.~K. Reid, J.~A. Scott, and K.~Turner}, {\em
  The factorization of sparse symmetric indefinite matrices}, IMA Journal of
  Numerical Analysis, 11 (1991), pp.~181--204.

\bibitem{ElBTZ96}
{\sc A.~{El-Bakry}, R.~A. Tapia, and Y.~Zhang}, {\em On convergence rate of
  newton interior-point algorithms in the absence of strict complementarity},
  Computational Optimization and Applications, 6 (1996), pp.~157--167.

\bibitem{ForGS94}
{\sc A.~Forsgren, P.~Gill, and J.~Shinnerl}, {\em Stability of symmetric
  ill-conditioned systems arising in interior methods for constrained
  optimization}, SIAM Journal on Matrix Analysis and Applications, 17 (1996),
  pp.~187--211.

\bibitem{FouM92}
{\sc R.~Fourer and S.~Mehrotra}, {\em Solving symmetric indefinite systems in
  an interior-point method for linear programming}, Mathematical Programming,
  62 (1993), pp.~15--39.

\bibitem{Gau77}
{\sc J.~Gauvin}, {\em A necessary and sufficient regularity condition to have
  bounded multipliers in nonconvex programming}, Mathematical Programming, 12
  (1977), pp.~136--138.

\bibitem{GolV96}
{\sc G.~H. Golub and C.~F. {Van Loan}}, {\em Matrix Computations}, The Johns
  Hopkins University Press, {B}altimore, third~ed., 1996.

\bibitem{Gou86}
{\sc N.~I.~M. Gould}, {\em On the accurate determination of search directions
  for simple differentiable penalty functions}, IMA Journal of Numerical
  Analysis, 6 (1986), pp.~357--372.

\bibitem{GouOST00}
{\sc N.~I.~M. Gould, D.~Orban, A.~Sartanaer, and P.~Toint}, {\em Superlinear
  convergence of primal-dual interior-point algorithms for nonlinear
  programming}, Technical Report TR/PA/00/20, CERFACS, April 2000.

\bibitem{Hag97a}
{\sc W.~W. Hager}, {\em Stabilized sequential quadratic programming},
  Computational Optimization and Applications, 12 (1999), pp.~253--273.

\bibitem{Hig96}
{\sc N.~J. Higham}, {\em Accuracy and Stability of Numerical Algorithms}, SIAM
  Publications, Philadelphia, 1996.

\bibitem{Hig97}
\leavevmode\vrule height 2pt depth -1.6pt width 23pt, {\em Stability of the
  diagonal pivoting method with partial pivoting}, SIAM Journal on Matrix
  Analysis and Applications, 18 (1997), pp.~52--65.

\bibitem{ManF67}
{\sc O.~L. Mangasarian and S.~Fromovitz}, {\em The {Fritz-John} necessary
  optimality conditions in the presence of equality and inequality
  constraints}, Journal of Mathematical Analysis and Applications, 17 (1967),
  pp.~37--47.

\bibitem{Meh92a}
{\sc S.~Mehrotra}, {\em On the implementation of a primal-dual interior point
  method}, SIAM Journal on Optimization, 2 (1992), pp.~575--601.

\bibitem{MonW93a}
{\sc R.~D.~C. Monteiro and S.~J. Wright}, {\em Local convergence of
  interior-point algorithms for degenerate monotone {LCP}}, Computational
  Optimization and Applications, 3 (1994), pp.~131--155.

\bibitem{RalW96}
{\sc D.~Ralph and S.~J. Wright}, {\em Superlinear convergence of an
  interior-point method for monotone variational inequalities}, in
  Complementarity and Variational Problems: State of the Art, M.~C. Ferris and
  J.~Pang, eds., SIAM Publications, Philadelphia, Penn., 1997, pp.~345--385.

\bibitem{RalW96b}
{\sc D.~Ralph and S.~J. Wright}, {\em Superlinear convergence of an
  interior-point method despite dependent constraints}, Mathematics of
  Operations Research, 25 (2000), pp.~179--194.

\bibitem{MWri98a}
{\sc M.~H. Wright}, {\em Ill-conditioning and computational error in interior
  methods for nonlinear programming}, SIAM Journal on Optimization, 9 (1998),
  pp.~84--111.

\bibitem{Wri93c}
{\sc S.~J. Wright}, {\em Stability of linear equations solvers in
  interior-point methods}, SIAM Journal on Matrix Analysis and Applications, 16
  (1994), pp.~1287--1307.

\bibitem{Wri97e}
\leavevmode\vrule height 2pt depth -1.6pt width 23pt, {\em Modifying {SQP} for
  degenerate problems}, Preprint ANL/MCS-P699-1097, Mathematics and Computer
  Science Division, Argonne National Laboratory, Argonne, Ill., 1997.
\newblock Revised June 2000.

\bibitem{IPPD96}
\leavevmode\vrule height 2pt depth -1.6pt width 23pt, {\em Primal-Dual
  Interior-Point Methods}, SIAM Publications, Philadelphia, 1997.

\bibitem{Wri94a_rev}
\leavevmode\vrule height 2pt depth -1.6pt width 23pt, {\em Stability of
  augmented system factorizations in interior-point methods}, SIAM Journal on
  Matrix Analysis and Applications, 18 (1997), pp.~191--222.

\bibitem{Wri98b}
\leavevmode\vrule height 2pt depth -1.6pt width 23pt, {\em Effects of
  finite-precision arithmetic on interior-point methods for nonlinear
  programming}, Preprint ANL/MCS-P705-0198, Mathematics and Computer Science
  Division, Argonne National Laboratory, Argonne, Ill., January 1998.

\bibitem{Wri98a}
\leavevmode\vrule height 2pt depth -1.6pt width 23pt, {\em Superlinear
  convergence of a stabilized {SQP} method to a degenerate solution},
  Computational Optimization and Applications, 11 (1998), pp.~253--275.

\bibitem{Wri96b}
\leavevmode\vrule height 2pt depth -1.6pt width 23pt, {\em Modified {C}holesky
  factorizations in interior-point algorithms for linear programming}, SIAM
  Journal on Optimization, 9 (1999), pp.~1159--1191.

\bibitem{SJW32}
{\sc S.~J. Wright and D.~Ralph}, {\em A superlinear infeasible-interior-point
  algorithm for monotone nonlinear complementarity problems}, Mathematics of
  Operations Research, 21 (1996), pp.~815--838.

\end{thebibliography}
\bibliographystyle{siam}

\end{document}